\documentclass[10pt,a4paper,twoside]{article}
\usepackage{color}
\usepackage{amsthm,amsmath,amssymb}
\usepackage{mathrsfs}

\definecolor{colorMH}{rgb}{0.,0.67,0}
\definecolor{colorBN}{rgb}{0.67,0.,0}
\definecolor{colorJV}{rgb}{0.,0.,0.67}
\definecolor{colorRH}{rgb}{0.,0.67,0.67}
\newtheorem{teorema}{Theorem}[section]
\newtheorem{lemma}[teorema]{Lemma}
\newtheorem{osse}[teorema]{Remark}
\newtheorem{prop}[teorema]{Proposition}

\newtheorem{coro}[teorema]{Corollary}
\newtheorem{prob}[teorema]{Problem}
\theoremstyle{definition}
\newtheorem{defin}[teorema]{Definition}
\newtheorem{remark}[teorema]{Remark}

\newcommand{\bele}{\begin{lemm}\begin{sl}}
\newcommand{\enle}{\end{sl}\end{lemm}}
\newcommand{\bedef}{\begin{defi}\begin{sl}}
\newcommand{\eddef}{\end{sl}\end{defi}}
\newcommand{\bete}{\begin{teor}\begin{sl}}
\newcommand{\ente}{\end{sl}\end{teor}}
\newcommand{\beos}{\begin{osse}\begin{rm}}
\newcommand{\eddos}{\end{rm}\end{osse}}
\newcommand{\bepr}{\begin{prop}\begin{sl}}
\newcommand{\empr}{\end{sl}\end{prop}}
\newcommand{\bepro}{\begin{prob}\begin{rm}}
\newcommand{\empro}{\end{rm}\end{prob}}
\newcommand{\bede}{\begin{defin}\begin{sl}}
\newcommand{\edde}{\end{sl}\end{defin}}
\newcommand{\beco}{\begin{coro}\begin{sl}}
\newcommand{\enco}{\end{sl}\end{coro}}
\newcommand{\disp}{\displaystyle}

\newcommand{\beeq}[1]{\begin{equation}\label{#1}}
\newcommand{\eddeq}{\end{equation}}

\newcommand{\beeqa}[1]{\begin{eqnarray}\label{#1}}
\newcommand{\eddeqa}{\end{eqnarray}}

\newcommand{\beal}[1]{\begin{align}\label{#1}}
\newcommand{\eddal}{\end{align}}

\newcommand{\bespl}[1]{\begin{split}\label{#1}}
\newcommand{\edspl}{\end{split}}

\newcommand{\bega}[1]{\begin{gather}\label{#1}}
\newcommand{\edga}{\end{gather}}

\newcommand{\beeqax}{\begin{eqnarray*}}
\newcommand{\eddeqax}{\end{eqnarray*}}

\def\qed{\ifmmode 
  \else \leavevmode\unskip\penalty9999 \hbox{}\nobreak\hfill
  \fi
  \quad\hbox{\hskip.5em\vrule width.4em height.6em depth.05em\hskip.1em}}
\def\endproofsym{\qed}
\renewenvironment{proof}[1][Proof]{\trivlist\item[\hskip\labelsep{\hskip0pt
    {\normalfont\scshape#1.}\hskip .321429\parindent}]\ignorespaces}
{\endproofsym\endtrivlist}
\def\endnobox{\def\endproofsym{}\end{proof}\def\endproofsym{\qed}}

\newcommand{\no}{\nonumber}

\newcommand{\beeqao}{\begin{eqnarray}\no}
\newcommand{\bealo}{\begin{align}\no}
\newcommand{\besplo}{\begin{split}\no}
\newcommand{\begao}{\begin{gather}\no}

\newcommand{\+}{\hspace{1pt}}
\newcommand{\perogni}{\forall\,}

\newcommand{\epsi}{\varepsilon}

\newcommand{\vf}{\varphi}

\newcommand{\cK}{\mathop{K}^{\circ}}

\DeclareMathOperator{\dist}{dist}

\DeclareMathOperator{\sgn}{sgn}

\def\({\left(}
\def\){\right)}
\def\Bbb{\mathbb}

\let\TeXchi\chi
\def\chi{{\setbox0 \hbox{\mathsurround0pt
$\TeXchi$}\hbox{\raise\dp0 \copy0 }}}

\newcommand{\teta}{\vartheta}

\newcommand{\K}{\mathscr{K}}

\def\R{\Bbb R}
\def\Dt{\partial_t}
\def\Dx{\Delta_x}
\def\Nx{\nabla_x}
\def\eb{\varepsilon}
\def\sgn{\operatorname{sgn}}
\def\dist{\operatorname{dist}}
\def\Cal{\mathcal}

\numberwithin{equation}{section}
\def\Dt{\partial_t}
\begin{document}

\title{\Large Finite-dimensional global and exponential attractors for the reaction-diffusion problem with an obstacle potential}
\author{Antonio Segatti\\
{\sl Dipartimento di Matematica ''F.Casorati''}\\
{\sl Universit\`a  di Pavia}\\
{\sl  Via Ferrata 1, 27100 Pavia, Italy}\\
{\sl antonio.segatti@unipv.it}\\
\and
Sergey Zelik\\
{\sl Departement of Mathematics}\\
{\sl University of Surrey}\\
{\sl Guildford, GU2 7XH, UK}\\
{\sl s.zelik@surrey.ac.uk}	
}
\date{}
\maketitle

\begin{abstract}
A reaction-diffusion problem with an obstacle potential
is considered in a bounded domain of $\R^N$. Under the assumption that the obstacle $\K$
is a closed convex and bounded subset of $\mathbb{R}^n$ with smooth boundary
or it is a closed $n$-dimensional simplex, we prove that the long-time behavior of
the solution semigroup associated with this problem can be described in terms of an exponential
attractor. In particular, the latter means that the fractal dimension of the associated global
attractor is also finite.
\end{abstract}

\section{Introduction}

This paper is devoted to the long-time behavior of solutions of the following reaction-diffusion system with an obstacle
potential in a bounded and regular
domain $\Omega\subset\R^N$
\begin{equation}\label{eq}
\begin{cases}
\partial_t u -\Dx u + \partial I_{\mathscr{K}}(u) -\lambda u\ni 0,\\
u\big|_{t=0}=u_0,\ \ u\big|_{\partial\Omega}=0.
\end{cases}
\end{equation}
Here $u=(u_1(t,x),\cdots,u_n(t,x))$ is an unknown vector-valued function,  $\Dx$ is a Laplacian with
respect to the variable $x$, $\lambda>0$ is a given constant,
$\K$ is  a given bounded closed convex
set in $\R^n$ containing zero and
$\partial I_{\K}$ stands for the subdifferential of its indicator function $I_\K$:
\begin{equation}
I_\K(u):=\begin{cases} 0,\ \ u\in \K,\\
\infty,\  u\notin \K.
\end{cases}
\end{equation}
We recall that the  subdifferential $\partial I_\K$
consists in the
set of vectors in $\R^n$ such that $y\in \partial I_\K(x)$ if and only
if $(y,x-w)_{\mathbb{R}^n}\ge 0$ for any $w\in \K$, where $(\cdot,\cdot)_{\mathbb{R}^n}$ is the scalar product
in $\mathbb{R}^n$. It is also well
known that $\partial I_\K$ turns out to be a multivalued maximal monotone operator
(see \cite{Br73}, pg.25; see also Section \ref{s1} below for the rigorous definitions).
\par
Equations and systems of the type \eqref{eq} appear quite often in the mathematical analysis
of phase transitions models
or in reaction diffusion processes with constraints.
In the first physical situation, \eqref{eq}
rules the evolution of the so called
order parameter, which
is an  $n$-dimensional vector in the case of multicomponent systems
(see, e.g.,
\cite{CoHo93}).
Moreover, since
the order parameter $u$ is usually related to the pointwise proportions
of the $n$ independent phases shown by system under study, it is physically
reasonable the fact that it attains values only in a bounded (convex)
subset of $\R^n$, usually an $n$-dimensional simplex
\begin{equation}\label{phasesimplex}
\K:= \left\{p=(p_1,\ldots,p_n)\in \mathbb{R}^n \mbox{ such that } \sum_{i=1}^{n}p_1\le 1,\ \ p_i\ge0,\ i=1,\cdots,n\right\}.
\end{equation}
 This is motivated by the requirement
that no void nor overlapping should appear between the phases.
In particular, equation \eqref{eq} with $n=2$ and $\lambda = 0$ appears in the Fr\'emond models of
 shape memory alloys
(see \cite{Fr} also for other models of phase change showing the ubiquity
of subdifferential operators in this framework)
where the two components $u_1$ and $u_2$ denote the pointwise proportions of the
two martensitic variants.
Finally, in the simpler scalar case, i.e. $n=1$ and $\K=[0,1]$ (for instance), equation \eqref{eq} is usually referred as
Allen Cahn equation with double obstacle.
\par
The mathematical analysis of equations of the type \eqref{eq} (and more general equations associated with
maximal monotone
operators in Hilbert spaces) has attracted the attention
of researchers for many years. In the particular case of equation \eqref{eq}, results
concerning existence, approximation
and long time behavior of solutions (e.g., in terms of global attractors) are known and by now classic (without any sake of
completeness and referring only to results in the Hilbert space framework, we quote
\cite{Br72}, \cite{Br73}, \cite{Ba76}, \cite{Te}, \cite{NoSaVe}).
\par
However, to the best of our knowledge, the
finite/infinite dimensionality
of the global attractor associated with the obstacle problems has not been yet  understood
(even in the simplest case of Allen-Cahn equation with double obstacle). Indeed, the classical
machinery for proving the finite-dimensionality (in terms of fractal or/and Hausdorff dimension) of
the global attractor (which perfectly works in many  cases of  dissipative systems generated
by non-linear PDEs with {\it regular} non-linearities, see \cite{BV,Te} and references therein)
is based on the so-called volume contraction arguments and requires the associated solution
semigroup to be (uniformly quasi-) {\it differentiable} with respect to the initial data at least
on the attractor.
\par
Unfortunately, this differentiability condition is usually {\it violated} if the underlying PDE has
singularities or/and degenerations and, in particular, it is clearly violated for the obstacle
problems like \eqref{eq}. Thus, the classical scheme is not applicable here and this makes the problem
much more difficult and interesting. In fact,
it has been recently shown that,
contrary to the usual regular case,
the singular/degenerate
dissipative systems can easily generate {\it infinite-dimensional} attractors even in bounded domains.
For instance,
the global attractor of the degenerate analogue of the real Ginzburg-Landau equation
\begin{equation}\label{0.pm}
\Dt u=\Dx(u^3)+u-u^3, \ \ u\big|_{\Omega}=0
\end{equation}
is infinite-dimensional for any bounded domain $\Omega$ of $\R^N$ (thanks to the degeneration at $u=0$), see \cite{EfZ3}.
On the other hand, in recent years
the finite-dimensionality of the global attractor result has been established for many important classes of
degenerate/singular dissipative systems including Cahn-Hilliard equations with logarithmic
potentials (see \cite{MZ1}), porous media equations (under some natural restrictions which exclude the
example of \eqref{0.pm}, see \cite{EfZ3}), doubly non-linear parabolic equations of
different types (see \cite{MZ2,MZ3,EfZ1}), etc.
In these papers, the finite dimensionality of the global attractor is typically
a consequence of the existence of a more refined object called {\it exponential attractor},
whose existence proof is often based on proper forms of the so called
squeezing/smoothing property for the differences of solutions.
\par
We remind that the concept of exponential attractor has been introduced in \cite{EFNT} in order to overcome
two major drawbacks of the global attractors: the slow (uncontrollable) rate of attraction and the
sensitivity to perturbations. Roughly speaking, an exponential attractor (which always contains the global one)
is a compact {\it finite-dimensional}
set in the phase space which attracts {\it exponentially} fast the images of all bounded sets as time tends to infinity
(see Section \ref{s2} for the rigorous definition). Thus it turns out that, in contrast to global attractors,
the exponential attractors are much more robust to perturbations (usually, H\"older continuous with respect to the
perturbation parameter). Moreover, the rate of convergence to the exponential attractor can be
 controlled in term of physical parameters of the system considered, see \cite{EFNT} and the more recent survey \cite{MZsur}
(and also references therein) for more details. Finally, the finite dimensionality of the global attractor
immediately follows from the finite dimensionality of the exponential attractor.
\par
The main aim of the present paper is to extend the exponential attractors theory to some classes of reaction-diffusion
problems with {\it obstacle}  potentials. Although the methods based on the proper squeezing/smoothing property for
the differences of solutions do not require the differentiability with respect to the initial data and, in principle,
can be applied also to the obstacle problem \eqref{eq}, this application is far from being straightforward
in our situation since, to this
end, one needs to produce estimates for the difference between the Lagrange multipliers
(namely the selection
of the subdifferential $\partial I_\K(u)$ which turns the differential inclusion \eqref{eq}
into an equation)
associated with two solutions $u_1(t)$
and $u_2(t)$.
This kind of estimates, which form the  Assumption $\mathscr{L}$ (see
\eqref{1.eq1} and \eqref{2.lag} for the rigorous formulation), roughly speaking
look as follows
\begin{equation}\label{0.lm}
\int_0^1\|\partial I_\K(u_1(t))-\partial I_\K(u_2(t))\|_{L^1(\Omega)}\,dt\le C\|u_1(0)-u_2(0)\|_{L^2}
\end{equation}
where $u_1$ and $u_2$ denote two solutions starting from the proper absorbing set and, with a little abuse of notation, we refer to $\partial I_\K$ as it were single valued.
Such kind of estimates,  to the best of our knowledge,
do not seem to be already known.\\
\noindent
This paper is organized as follows
In Section \ref{s1}, we recall the basic results related with the wellposedness and regularity of solutions of the obstacle
problem \eqref{eq}. As usual, the results on the singular problem \eqref{eq} are obtained by approximating the singular
obstacle potential by more regular ones.
The usual regularization used for this kind of problems is the Moreau-Yosida approximation (see, e.g., \cite{Br73}).
However, it turns out that, in order to prove estimates of the type \eqref{0.lm}, it is more convenient to implement
different kind of approximation schemes. Thus, in Section \ref{s1} we also give a sketch of the wellposedness result for \eqref{eq}
(which is well known, see \cite{Br73} and \cite{Ba76})
by introducing another approximation scheme. This scheme becomes also very useful
to prove an $L^{\infty}$-estimate for the approximation of $\partial I_\K$ (independent of the approximation parameter) via the maximum pronciple. This
kind of estimate at the $\epsi$-level will guarantee that the same bound remains valid also for $\partial I_\K$.
This fact will be rather crucial in proving \eqref{0.lm}.
\par
Under the assumption that the
estimate \eqref{0.lm} for the Lagrange multipliers is known, in Section  \ref{s2}
we give the {\it conditional} proof of the existence of an exponential attractor
for the solution semigroup $S(t)$ associated with equation \eqref{eq}.
This result is obtained using some modification of
the so-called method of $l$-trajectories which was originally introduced by M\'alek and Ne{\v{c}}as
in \cite{MN} and is widely used nowadays in the  attractors theory, see \cite{MP,MZsur} and references therein.
\par
Section \ref{s3} is the key part of the paper and it is devoted to verifying estimate \eqref{0.lm}
for the case in which the boundary $\partial \K$ of the convex set $\K$ is regular enough.
The proof is based on the maximum principle and
the Kato inequality together with rather delicate construction
of the approximating potentials. Unfortunately, the most relevant (from the applications point of view)
choice of $\K$ is that of the simplex \eqref{phasesimplex} which does not fit with the regularity
assumptions on the boundary of $\K$. Thus, we have to consider this particular case separately using a
{\it different} approximation scheme for the subdifferential $\partial I_\K$.
As a result, the existence of an exponential attractor is proved for
the case in which the boundary of $\K$ is smooth and for the particular choice of $\K$ as a simplex \eqref{phasesimplex}
(although we expect that this result should be true for any closed and bounded convex set $\K$).
The section closes with some discussion
on possible extensions of our result to more general classes of reaction-diffusion equations.
\par
Finally, in the last Section \ref{s4}, we study the convergence of the exponential attractors $\Cal M_\eb$ for the
regular approximating problems to the limit exponential attractor $\Cal M_0$ of the singular problem \eqref{eq}.
For simplicity, we restrict ourselves there only to the case of the simplex \eqref{phasesimplex} and verify that
$$
\dist_{(L^\infty(\Omega))^n}^{sym}(\Cal M_\eb,\Cal M_0)\le C\eb^\kappa,
$$
where $\dist^{sym}$ stands for the symmetric Hausdorff distance and the positive constants $C$ and $\kappa$ are
independent of $\eb$. The result is based on the proper application of the abstract theorem on perturbations of
 exponential attractors proved in \cite{FGMZ} to the obstacle problem \eqref{eq}.

\section{Well-posedness and regularity}\label{s1}
The aim of this section is to recall some known facts about the solutions of
reaction-diffusion problem \eqref{eq} with obstacle potential and to formulate some additional estimates which will
be crucial for what follows. We start with
\\
\noindent \textbf{A word on the notation:} The unknown function $u$ is actually a vector
valued function but, for the sake of simplicity, will be denoted as a scalar valued function.
Consequently, also the functional spaces we will use in the course of the paper we will
have a "scalar" notation. This means that a notation like, e.g., $L^2$
will be preferred to a (more precise) notation like $(L^2(\Omega))^n$. The same applies to dualities ($\langle\cdot,\cdot\rangle$)
and scalar products ($(\cdot,\cdot)$). Moreover, we will indicate with same
symbols $\K$ and $I_\K$ the convex in $\mathbb{R}^n$ and its indicator function and their realization in $L^2$ (see the next proposition \ref{Prop1.point}).
Thus,
the definition of (weak) solutions of our problem is.
\begin{defin}\label{Def1.sol} A function $u=u(t,x)$ is a solution of the obstacle problem
\eqref{eq} if  $u(t,x)\in \K$ for almost all $(t,x)\in[0,T]\times\Omega$,
\begin{equation}\label{1.sp}
u\in C([0,T];L^2)\cap L^2(0,T;H^{1}_{0}), \ \
\Dt u\in L^2(0,T;H^{-1}),
\end{equation}
and the following variational inequality holds for almost every $t\in(0,T]$
\begin{equation}\label{1.vareq}
 \langle \Dt u(t), u(t)-z\rangle  + (\nabla u,\nabla (u-z)) \le \lambda(u,u-z), \hbox{ for any } z\in H^1_{0}\cap \K.
\end{equation}
\end{defin}
The next theorem is a standard result in the theory of the evolution equations associated
with maximal monotone operators (see the seminal references \cite{Br72}, \cite{Br73} and \cite{Ba76}).
\begin{teorema}\label{Th1.wp}[Well posedness]
Let $\K$ be a closed and bounded convex set containing the origin in $\mathbb{R}^n$. The,
for any given measurable $u_0$ such that $u_0(x)\in \K$ for almost all $x\in\Omega$
there exists a unique global solution $u$ of problem \eqref{eq} in the sense
of definition \ref{Def1.sol}. Moreover,
$u(t)\in H^2$ and $\Dt u(t)\in L^2$ for $t>0$
and the following estimate holds:
\begin{equation}\label{1.reg}
\|\Dt u(t)\|_{L^2}+\|u(t)\|_{H^2}\le C\frac {t+1}t,
\end{equation}
where the constant $C$ depends on $\lambda$ and $\K$, but is independent of $u_0$. If, in addition,
$u_0\in H^2\cap H^1_{0}$,
then the term with $\frac{1+t}{t}$ in the right-hand side of \eqref{1.reg}
can be removed.
\end{teorema}
\begin{proof}
 Although this result is well-known, for the convenience of the reader,
we briefly recall one of its possible proof and deduce the regularity estimate \eqref{1.reg}.
\par
We first note that the uniqueness follows immediately from the variational
inequality \eqref{1.vareq}. Indeed, let $u$ and $v$ be two solutions of \eqref{eq}.
Then, taking a sum of \eqref{1.vareq} for $u(t)$ and $z=v(t)$ with \eqref{1.vareq} for
$v(t)$ with $z=u(t)$, we have
$$
(\Dt u(t)-\Dt v(t),u(t)-v(t))+\|\Nx(u(t)-v(t))\|^2_{L^2}\le\lambda\|u(t)-v(t)\|^2_{L^2}.
$$
Denoting now $w(t):=u(t)-v(t)$ we arrive at the differential inequality
\begin{equation}\label{1.gr}
\frac{1}{2}\frac d{dt}\|w(t)\|^2_{L^2}+\|\Nx w(t)\|^2_{L^2}\le\lambda\|w(t)\|^2_{L^2}.
\end{equation}
Applying the Gronwall inequality and using that, by definition $u, v\in C^{0}([0,T];L^2)$
for any $T>0$, we finally have
\begin{equation}\label{1.lip}
\|w(t)\|^2_{L^2}+\int_0^t\|w(s)\|^2_{H^1_{0}}\,ds\le Ce^{\lambda s}\|w(0)\|^2_{L^2},
\end{equation}
where the constant $C$ is independent of $u$ and $v$. Thus, the uniqueness holds.
\par
The existence part is slightly more delicate and requires the approximations of
the singular convex potential $I_\K$ by suitable regular ones $F_\eb$.
In fact, thanks to the uniqueness result this approximation can be done in several ways.
The typical choice is usually the Moreau-Yosida approximation (see \cite{Br73}).
However, in view of the next results (see in particular Proposition \ref{Prop1.linf} and Theorems \ref{Th3.smooth} and \ref{Th3.simplex}),
it is more convenient to adopt another kind of approximation.
To this end, we let $M(u)$ be the distance from
the point $u\in\R^n$ to the convex set $\K$, namely the real valued function
$M$
\begin{equation}\label{1.d}
M(u):=\operatorname{dist}(u,\K).
\end{equation}
Then, $M$ is convex, globally Lipschitz continuous
and smooth outside $\K$. In addition there holds that
\begin{subequations}\label{1.dd}
\begin{align}
&\disp M(u)\ge0,\ M(u)=0, \hbox{ if } u\in \K,\label{1.dd1}\\
&\disp  |\nabla M(u)|=1, \hbox{ if } u\notin \K.\label{1.dd2}
\end{align}
\end{subequations}
Now for any $\eb>0$, we introduce the real function
\begin{equation}\label{1.ddd}
f_\eb(z):=\begin{cases} 0,\ \ z\le\eb,\\
                        \eb^{-1}(z-\eb)^2,\ \ z\ge\eb.
          \end{cases}
\end{equation}
Finally, the desired approximation is defined a
\begin{equation}\label{1.dddd}
F_\eb(u):=f_\eb(M(u)).
\end{equation}
Then, obviously, $F_\eb(u)$ is convex and smooth (at least $C^{1,1}$) and $F_\eb(u)=F'_\eb(u) =0$
for all $u\in \K$, where $F_\epsi'$ denotes the gradient of $F_\epsi$. We thus consider
the following approximation of the problem \eqref{eq}:
\begin{equation}\label{1.eqreg}
\begin{cases}
\Dt u-\Dx u+F'_\eb(u)=\lambda u,\\
u\big|_{\partial\Omega}=0,\ \ u\big|_{t=0}=u_0,
\end{cases}
\end{equation}
where the initial data $u_0\in L^\infty(\Omega)$ and $u_0(x)\in \K$ a.e. in $\Omega$.
\par
Our next step is to obtain a number of uniform (with respect to $\eb\to0$) estimates
for the solutions of the regular problems \eqref{1.eqreg}. We start with the usual
$L^2$-estimate. Indeed, multiplying equation \eqref{1.eqreg} by $u$, integrating
over $\Omega$ and using the obvious
facts that
$$
\frac{1}{2}(F'_\eb(u),u)-\lambda \|u\|^2_{L^2}\ge -C,\ \  (F'_\eb(u),u)=(F'_\eb(u)-F'_\eb(0),u)\ge0
$$
where the constant $C$ depends only on $\lambda$ and on $\Omega$ but is independent of $\eb\searrow 0$,
we arrive at
$$
\frac{1}{2}\frac d{dt}\|u(t)\|^2_{L^2}+\|\Nx u(t)\|^2_{L^2}+\frac{1}{2}(F'_\eb(u),u)\le C.
$$
Applying the Gronwall inequality to this relation, we see that
\begin{equation}\label{1.l2}
\|u(t)\|^2_{L^2}+\int_t^{t+1}\|u(s)\|^2_{H^1_{0}}+(F'_\eb(u(s)),u(s))\,ds\le C.
\end{equation}
Here we have implicitly used that $\K$ is bounded and, therefore, $\|u_0\|_{L^2}$
is uniformly bounded for all admissible initial data $u_0$.
\par
Next, we obtain the uniform energy estimate for the solutions of \eqref{1.eqreg}.
To this end, we multiply the equation \eqref{1.eqreg} by $\Dt u(t)$ and
 integrate over $x$:
$$
\frac {d}{dt}\Big(\|\Nx u(t)\|^2_{L^2}+2(F_\eb(u(t)),1)\Big)
+\|\Dt u(t)\|^2_{L^2}\le\lambda \|u(t)\|^2_{L^2}.
$$
Integrating this estimate with respect to time, using \eqref{1.l2} together with the obvious estimate
$$
0\le (F_\eb(u),1)\le (F'_\eb(u),u)
$$
(which is an immediate consequence of the convexity of $F_\eb$) and arguing in a standard way
(see, e.g., \cite{Z}), we deduce the desired $H^1$-energy estimate
\begin{equation}\label{1.h1}
\|u(t)\|^2_{H^1_{0}}+(F_\eb(u(t)),1)+\int_{t}^{t+1}\|\Dt u(s)\|^2_{L^2}\,ds\le C\frac{t+1}t,\ \ \ t>0,
\end{equation}
where the constant $C$ is independent of $u_0$ and $\eb$.
\par
We are now ready to verify the analogue of \eqref{1.reg} for the approximate solutions.
We start with the estimate of $\Dt u(t)$. To this end, we differentiate \eqref{1.eqreg}
with respect to time and denote $v(t):=\Dt u(t)$. Then, we have
$$
\Dt v-\Dx v+F_\eb''(u)v=\lambda\Dt u.
$$
Multiplying this equation by $v(t)$ integrating with respect to time, using \eqref{1.h1} together with the monotonicity of $F'_\eb(u)$
and arguing in a standard way (see again \cite{Z} for details), we have
\begin{equation}\label{1.dt}
\|\Dt u(t)\|^2_{L^2}\le C\frac{(1+t)^2}{t^2}
\end{equation}
for some constant $C$ independent of $u_0$ and $\eb$.
\par
Now, in order to obtain the $H^2$-part of estimate \eqref{1.reg}, it is sufficient
to rewrite equation \eqref{1.eqreg} in the form of elliptic equation (for every fixed $t$)
 \begin{equation}\label{1.ell}
\Dx u(t)-F_\eb'(u(t))=\Dt u(t)-\lambda u(t),\ \ u(t)\big|_{\partial\Omega}=0,
\end{equation}
multiply it by $\Dx u(t)$, integrate over $\Omega$
and use the fact that, thanks to the convexity of $F_\epsi$, the term $(F'_\epsi(u),-\Delta_x u)$ is nonnegative.
Then, due to \eqref{1.dt}, \eqref{1.h1} and
the elliptic regularity result for the Laplacian (recall that the domain $\Omega$ is assumed to be smooth enough),
\begin{equation}\label{1.h2}
\|u(t)\|_{H^2}\le C\|\Dx u(t)\|_{L^2}\le C(\|\Dt u(t)\|_{L^2}+\lambda\|u(t)\|_{L^2})\le
C'\frac{t+1}{t},
\end{equation}
where the constants $C'$ and $C$ are independent of $\eb$ and $u_0$. Furthermore,
from equation \eqref{1.ell}, we now conclude that
$$
\|F'_\eb(u(t))\|_{L^2}\le C\frac{t+1}t.
$$
Finally, we are now ready to pass to the limit $\eb\to0$ in equations \eqref{1.eqreg}
and verify the existence of a solution $u$ for the limit obstacle problem. Indeed,
let $\eb_n>0$ be a sequence of positive numbers such that $\eb_n\to0$ as $n\to\infty$.
Assume for the first that the initial data $u_0$ is smooth enough and verify the constraint, namely
$$
u_0\in H^2\cap H^1_0 \hbox{ and }  u_0(x)\in \K,\ \ \hbox{a.e. in } \in\Omega.
$$
Let $u_n(t)=u_{\eb_n}(t)$ be the corresponding solution of the approximating problem
\eqref{1.eqreg}.
Then, analogously to the derivation of \eqref{1.reg} for $u_n$, but using in addition
that $u_0$ is regular, we have
\begin{equation}\label{1.unest}
\|\Dt u_n(t)\|_{L^2}+\|u_n(t)\|_{H^2}+\|F_{\eb_n}'(u_n(t))\|_{L^2}+\|F_\eb(u_n(t))\|_{L^1}\le Q(\|u_0\|_{H^2}),
\end{equation}
 where the function $Q$ is independent of $n$ and $t\ge0$. Thus, up to not relabeled subsequence, we have
 that $u_n$ converges weakly star in
 $$
 X:=L^\infty(0,T;H^2)\cap W^{1,\infty}(0,T;L^2)
 $$
 and strongly (thanks to \cite[Cor. 4]{simon}) in $C^0([0,T]; H^1_0)$
 to some function $u\in X \cap C^0([0,T]; H^1_0)$.
 So, we only need to check that $u$ is a desired solution
 of \eqref{eq} in the sense of Definition \ref{Def1.sol}. Indeed, the regularity of
 $u$ is obvious. The fact that $u(t,x)\in \K$ for almost $(t,x)\in[0,T]\times\Omega$
 follows in a standard way from the uniform bounds for $F_{\eb_n}(u_n(t))$ in the $L^1$-norm
 and from the fact that $F_{\eb_n}(w)\to\infty$ for all $w\notin \K$.
So, we only need to check the variational inequality \eqref{1.vareq}. Let $z=z(x)$ be any
admissible test function. Then, using the monotonicity of $F'_\eb(u)$ and the fact
that $F'_\eb(v)=0$ for all $v\in \K$, we see that
$$
(F'_\eb(u(t)),u(t)-z)=(F'_\eb(u(t))-F'_\eb(z),u-z)\ge0.
$$
Multiplying now \eqref{1.eqreg} by $u(t)-z$, integrating over $\Omega$ and using the
last inequality, we arrive at
$$
\langle \Dt u_n(t),u_n(t)-z\rangle +(\Nx u_n(t),\Nx(u_n(t)-z))\le\lambda(u_n(t),u_n(t)-z).
$$
Passing to the limit $n\to\infty$ in this inequality and using the above convergences,
we see that \eqref{1.vareq} holds
for almost any $t$ and all admissible test functions $z$. Thus, we have proved the
existence of the desired solution $u$ for the case of smooth initial data $u_0$.
\par
For general initial data $u_0$, the corresponding solution can be now constructed by
approximating the non-smooth initial data $u_0$ by smooth ones $u_0^n$ and passing to
the limit $n\to\infty$. To conclude, we finally recall that the continuity
of the limit function near zero ($u\in C([0,T],L^2)$) follows from the global Lipschitz continuity \eqref{1.lip}).
Theorem \ref{Th1.wp} is proved.
\end{proof}
\begin{remark}\label{Rem1.app} The above proof reveals that the concrete form of the
functionals $F_\eb(u)$ is not essential for the proof of the existence result.
In fact, the solutions $u_\eb$ of the approximating problems \eqref{1.eqreg}
will converge to the unique solution $u$ of the obstacle problem if the following
assumptions are satisfied:
\begin{enumerate}
\item The functions $F_\eb $ are {\it convex} and regular enough\\
\item  $F'_\eb \to \partial I_\K$ as $\epsi \searrow 0$.
\end{enumerate}
Where the last condition means that
\begin{eqnarray*}
& i)& \;\;F'_\eb(v)\to 0,\ \ \forall v\in \K,\\
  & ii)& \;\; \vert F'_\eb(v)\vert \to+\infty,\ \ \forall v\notin \K.
\end{eqnarray*}
In other words, we only need the subdifferential $F'_{\epsi}(u)$
of $F_\epsi$ to be an approximation, via graph convergence, of the subdifferential $\partial I_\K$
(see, e.g., \cite[Proposition 3.60 and Theorem 3.66]{attouch}).
However, as we have already mentioned, the choice of suitable approximating functionals $F_\eb(u)$ turns out to be crucial
for our method of estimating the dimension of the global attractor. Thus we will construct this approximation
in a rather specific way depending on the structure of the convex set $\K$ (see, however, the next sections).
\end{remark}
For the long time analysis it is actually more convenient
to reformulate the variational inequality \eqref{1.vareq} in terms of an equation
coupled with a differential inclusion for the subdifferential of the indicator function of $\K$ (actually of its realization). This
subdifferential, which a priori should be undestood with respect to the duality
$H^1_0$-$H^{-1}$, actually makes sense in $L^2$ and thus almost everywhere in $\Omega$ (see Theorem
\ref{Th1.wp} and the next Proposition \ref{Prop1.point}).
Thus, we introduce the function (named Lagrange multiplier in what follows)
\begin{equation}\label{1.h}
h_u(t):=-\Dt u(t)+\Dx u(t)+\lambda u, \;\;\; h_u(t) \in \partial I_\K(u).
\end{equation}
Then, due to Theorem \ref{Th1.wp}, $h_u\in L^\infty(\tau,T;L^2)$, for any $\tau>0$. Moreover, the definition of subdifferential (w.r.t. the
$L^2$ scalar product) gives
\eqref{1.vareq},
\begin{equation}\label{1.vein}
(h_u(t),u(t)-z)\ge0,\ \ \text{for almost all $t\in[0,T]$}
\end{equation}
and any admissible test function $z=z(x)$. The last inequality can be also written in a point-wise form.
\begin{prop}\label{Prop1.point} Let the assumptions of Theorem \ref{Th1.wp} hold and let $h_u=h_u(t,x)$
be the Lagrange multiplier associated with the solution $u(t)$ of problem \eqref{eq}. Then
\begin{equation}\label{1.varep}
(h_u(t,x),u(t,x)-Z)_{\mathbb{R}^{n}}\ge 0,\ \ \forall Z\in \K
\end{equation}
and almost all $(t,x)\in[0,T]\times\Omega$, which means
$$
h_u(t,x)\in\partial I_\K(u(t,x)),\ \ \text{ a.e. in}\  [0,T]\times\Omega.
$$
\end{prop}
\begin{proof} Indeed, arguing as in the proof of Theorem \ref{Th1.wp}, we see that
$$
h_{u_n}(t):=F'_{\eb_n}(u_n(t))\to h_u(t)
$$
weakly-star in the space $L^\infty(\tau,T; L^2)$, for any $\tau>0$. Thus, it is sufficient to verify \eqref{1.varep}
for the functions $h_{u_n}$ only. But these inequalities are immediate due to the monotonicity of $F'_\eb$ and
the fact that $F'_\eb(Z)=0$ if $Z\in \K$. Indeed,
$$
(h_{u_n}(t,x),u_n(t,x)-Z)_{\mathbb{R}^n}=(F'_{\eb_n}(u_n(t,x))-F'_{\eb_n}(Z)),u_n(t,x)-Z)_{\mathbb{R}^n}\ge0.
$$
Passing to the limit $n\to\infty$ in that inequalities, we deduce \eqref{1.varep} and finish the proof of the proposition.\
\end{proof}
Thus, the obstacle problem \eqref{eq} can be rewritten in terms of functions $u$ and $h_u$ as follows:
\begin{equation}\label{1.eq1}
\begin{cases}
\Dt u-\Dx u+h_u=\lambda u, \ \ \text{in the sense of distributions,}\\
h_u(t,x)\in \partial I_\K(u(t,x)),\ \ \text{for almost all $(t,x)\in\R_+\times\Omega$},\\
u\big|_{t=0}=u_0,\ \ u\big|_{\partial\Omega}=0.
\end{cases}
\end{equation}
The next proposition shows that the function $h_u$ is, in a fact, globally bounded in the $L^\infty$-norm.

\begin{prop}\label{Prop1.linf} Let the assumptions of Theorem \ref{Th1.wp} hold and let $h_u(t)$ be the Lagrange
multiplier associated with the solution $u(t)$ of problem \eqref{eq}. Then, $h_u\in L^\infty(\R_+\times\Omega)$ and
\begin{equation}\label{1.linf}
\|h_u(t)\|_{L^\infty}\le C,
\end{equation}
where the constant $C$ depends only on $\K$ and $\lambda$ (and is independent of $u$ and $t\ge0$).
\end{prop}
\begin{proof}
Thanks to the lower semicontinuity of norms with respect to the weak convergence,
we will verify \eqref{1.linf} only in the case in which $h_u$  in equation \eqref{1.eq1} is replaced
by its approximation  $h_{u_n}(t):=F'_{\eb_n}(u_n(t))$.
To this end, we test equation \eqref{1.eqreg} in the scalar product of $\mathbb{R}^n$
with $\nabla M(u)$ (where the function $M$ is the same as in the proof of Theorem \ref{Th1.wp}),
and use that, due to the convexity of
$M$,
$$
(-\Dx u, \nabla M(u))_{\mathbb{R}^n}=-\Dx(M(u)) + (\mathcal{H}(M(u))\Nx u,\Nx u)_{\mathbb{R}^n}\ge - \Dx(M(u))
$$
where $\mathcal{H}(M)$ denotes the Hessian matrix of $M$
(actually, $\mathcal{H}(M)(u)$ does not exist if $u\in\partial \K$, but the {\it inequality} still holds and can be easily
verified, say, by approximating the non-smooth convex function $M$ by the smooth convex ones).
Then, we arrive at the
differential inequality for the function $M(u)$
$$
\Dt M(u)-\Dx M(u) + f'_{\eb_n}(M(u))\le \lambda (u,\nabla M(u)),\ \ M(u)\big|_{t=0}=0,
$$
where we have implicitly used that $|\nabla M(u)|^2=1$ for $ u\notin \K$ (see \eqref{1.dd2}). Furthermore, since $|u|\le M(u)+C_\K$, where
$C_\K:=\operatorname{diam}(\K)$, we finally have
\begin{equation}\label{1.max}
\Dt M(u)-\Dx M(u) + f'_{\eb_n}(M(u))\le \lambda M(u)+\lambda C_K,\ \ M(u)\big|_{t=0}=0.
\end{equation}
Applying the comparison principle to the scalar parabolic equation \eqref{1.max}, we see that
$$
M(u(t,x))\le v_{\eb_n},\ \ \hbox{ for a.e. } (t,x)\in\R_+\times\Omega,
$$
where the constant (w.r.t. $x$ and $t$) $v_{\eb_n}>0$ is the solution of the following equation
$$
f'_{\eb_n}(v_{\eb_n})=\lambda v_{\eb_n}+\lambda C_\K.
$$
In addition, from  \eqref{1.ddd}, we see that $v_\eb\to0$ as $\eb_n\to0$ (i.e. $n\nearrow +\infty$).
For this reason
$$
f'_{\eb_n}(M(u))\le f'_{\eb_n}(v_{\eb_n})=\lambda v_{\eb_n}+\lambda C_\K\le C_\lambda.
$$
and, therefore,
$$
|h_{u_n}(t)|=f'_{\eb_n}(M(u_n(t))\le C_\lambda
$$
uniformly with respect to $n\nearrow\infty$. Finally, passing to the limit $n\nearrow\infty$,
we arrive at \eqref{1.linf} (with $C=\lambda\operatorname{diam}(\K)$) and
finish the proof of the proposition.
\end{proof}
\begin{remark} In contrast to Theorem \ref{Th1.wp} and Proposition \ref{Prop1.point} which
are based only on the energy type estimates which are valid for much more
general equations, e.g., with non-scalar diffusion matrix, etc., the $L^\infty$-estimate
obtained in Proposition \ref{Prop1.linf} is based on the maximum/comparison principle and requires
the diffusion matrix to be scalar. In particular, we do not know whether or not this estimate remains
true even for the case of diagonal, but non-scalar diffusion matrix.
\end{remark}
As direct consequence of the previous Proposition, we have that \eqref{1.eq1} could
be understood as the heat equation
$$
\Dt u-\Dx u=\lambda u-h_u
$$
with the external forces belonging to $L^\infty$. Thus, the parabolic interior regularity
estimates give (see. e.g.,  \cite{LSU})
\begin{coro}\label{Cor1.c2} Let the assumptions of Theorem \eqref{Th1.wp} hold and let $u(t)$ be a solution of problem \eqref{eq}.
Then, for every $\nu>0$, $u(t)\in C^{2-\nu}(\Omega)$ for $t>0$ and the following estimate holds:
\begin{equation}\label{1.c2}
\|u(t)\|_{C^{2-\nu}}\le C_\nu\frac{1+t^\alpha}{t^\alpha},
\end{equation}
where the positive constants $C_\nu$ and $\alpha$ are independent of $t$ and $u$. If, in addition,
$u_0\in C^{2-\nu}\cap H^1_0$, then \eqref{1.c2} holds with $\alpha=0$.
\end{coro}


\section{Global and exponential attractors}\label{s2}
The aim of this section is to study the long-time behavior of solutions of problem \eqref{eq} in terms of global and exponential attractors. We first recall that, due to
Theorem \ref{Th1.wp}, equation \eqref{eq} generate a (dissipative) semigroup $\{S(t),\,t\ge0\}$ in the phase space
\begin{equation}\label{2.phase}
\Phi=\Phi_\K:=\{u\in L^\infty:\ \ u(x)\in \K \ \text{ for almost all $x\in\Omega$}\},
\end{equation}
i.e.,
\begin{equation}
S(t):\Phi\to \Phi, \ \ S(t)u_0:=u(t),
\end{equation}
where $u(t)$ is the solution to \eqref{eq} at
 time $t$. Moreover, due to estimate \eqref{1.lip} this semigroup is globally Lipschitz
 continuous in the $L^2$-metric
 \begin{equation}\label{2.lip}
 \|u_1(t)-u_2(t)\|_{L^2}^2+\int_0^t\|u_1(s)-u_2(s)\|_{H^1}^2\,ds\le
 Ce^{\mu t}\|u_0^{1}-u_0^{2}\|_{L^2}^2,
 \end{equation}
 where the positive constants $C$ and $\mu$ are independent of $t$ and of the initial data $u_0^{1},u_0^{2} \in \Phi$.
 In addition, due to Corollary \ref{Cor1.c2}, we have the following regularization estimate
 \begin{equation}\label{2.reg}
 \|\Dt u(t)\|_{L^2}+\|u(t)\|_{C^{2-\nu}(\Omega)}\le C_\nu\frac{1+t^\alpha}{t^\alpha},
 \end{equation}
 where $\nu>0$ is arbitray and the positive constants $C_\nu$ and $M$ are independent of
 the initial condition $u_0$ and of $t$. These two estimates immediately imply the existence of a global attractor
 $\Cal A$ for
 the semigroup $S(t)$ associated with the obstacle problem \eqref{eq}. We recall that, by
 definition, a set $\Cal A\subset\Phi$ is a global attractor for the semigroup $S(t):\Phi\to\Phi$ if
 \begin{enumerate}
 \item The set $\Cal A$ is compact in $\Phi$;
 \item  It is strictly invariant: $S(t)\Cal A=\Cal A$, $t\ge0$;
 \item For every neighborhood $\Cal O=\Cal O(\Cal A)$ of $\Cal A$ in $\Phi$ there exists a time
 $T=T(\Cal O)$ such that
 \end{enumerate}
 \begin{equation}\label{2.attr}
 S(t)\Phi\subset \Cal O(\Cal A),\ \ t\ge T.
 \end{equation}
\begin{remark}\label{Rem2.bound}
The attraction property \eqref{2.attr} is usually formulated not for the whole phase space
$\Phi$, but for the {\it bounded} subsets of it only. In our case, however, the whole phase space $\Phi$ is automatically bounded (since $\K$ is bounded in $\R^n$), so we need not to use bounded sets to define the attractor. Note also that the attraction property \eqref{2.attr} can be reformulated as follows:
\begin{equation}\label{2.dist}
\dist_{L^\infty}(S(t)\Phi,\Cal A)\to0\ \ \text{as $t\to\infty$},
\end{equation}
where $\disp \dist_V(X,Y):=\sup_{x\in X}\inf_{y\in Y}d_V(x,y)$ is the non-symmetric Hausdorff distance between the subsets $X$ and $Y$ of
the metric space $V$.
\end{remark}

\begin{teorema}\label{Th2.globattr}[Global Attractor]
Under the assumptions of Theorem \ref{Th1.wp}, the semigroup $S(t)$ associated with the obstacle equation \eqref{eq} possesses
the global attractor
$\mathcal{A}$ in $\Phi$ which is bounded in $C^{2-\nu}(\Omega)$ for every $\nu>0$. This
attractor is generated by all the trajectories of the semigroup $S(t)$ defined for all $t\in\R$:
\begin{equation}\label{2.str}
\Cal A={\rm K}\big|_{t=0},
\end{equation}
where ${\rm K}\subset L^\infty(\R,\Phi)$ is a set of all solutions of \eqref{eq} defined
for all $t\in\R$.
\end{teorema}
\begin{proof}
The proof of this theorem is standard and well known (see, for instance \cite{BV}, \cite{Te},
 or \cite{Lau94}). Indeed, due to estimate \eqref{2.reg} the semigroup $S(t)$ possesses
 an absorbing set which is bounded in $C^{2-\nu}(\Omega)$ and, therefore, compact in $L^\infty(\Omega)$
 and estimate \eqref{2.lip} guarantees that the semigroup has a closed
 graph. Thus, all assertions of the theorem follow from the abstract attractor existence
 criterium, see \cite{BV}.
\end{proof}
Recall that the semigroup $S(t)$ associated with
the obstacle problem \eqref{eq} possesses a global Lyapunov function of the form
\begin{equation}\label{2.lyap}
\Cal L(u):=\|\Nx u\|^2_{L^2}-\lambda\|u\|^2_{L^2}.
\end{equation}
Indeed, using the test function $z=u(t-h)$ in the variational inequality \eqref{1.vareq},
dividing it by $h>0$ and passing to the limit $h\to0$, we arrive at
$$
\|\Dt u(t)\|_{L^2}^2+\frac d{dt}\Cal L(u(t))\le 0.
$$
Therefore, according to the general theory (see e.g., \cite{BV}), every trajectory
$u(t)=S(t)u_0$ tends as $t\to\infty$ to the set $\Cal R$ of all equilibria of problem \eqref{eq}
$$
\dist(S(t)u_0,\Cal R)\to0,\ \text{ as $t\to\infty$}.
$$
However, in contrast to the case of regular systems, the equilibria set $\Cal R$
is generically {\it not discrete} for the  obstacle type singular problems.
Thus, in our situation, the existence of a Lyapunov function does not allow
to obtain the stabilization of every trajectory to a single equilibrium even
in "generic" situation. In addition, the semigroup $S(t)$ is  not differentiable
with respect to the initial data (it is in fact only globally Lipschitz continuous), so we are not able
to construct the stable/unstable manifolds associated with an equilibrium. Thus, the
so-called theory of {\it regular} attractors is not applicable to equations of the type \eqref{eq}.
Moreover, due to the above mentioned non-differentiability, the standard way of
proving the finite-dimensionality of the global attractor based on the volume contraction
method does not work here. So, the existence of the finite-dimensional reduction
for  the associated  long time dynamics becomes a non-trivial problem which,
to the best of our knowledge, it has not been yet tackled.
In this paper we will prove that the global attractor for \eqref{eq} has finite fractal dimension
by using the concept of the so-called {\it exponential}
attractor and the estimation of the dimension based on the proper chosen squeezing/smoothing property
for the difference between two solutions. This method has the advantage that it does not
require the differentiability with respect to the initial data
The existence of an exponential attractor is interest in itself. In fact, we
recall once more that the global attractor
represents the first (although extremely important)
step in the understanding of the long-time dynamics of a given
evolutive process.
However, it may also present some severe drawbacks.
Indeed, as simple examples
show,  the rate of convergence
to the global attractor may be arbitrarily slow. This fact makes the global attractor
very sensitive to perturbations and to numerical approximation.
In addition,  it is usually  extremely difficult
to estimate the rate of convergence to the global attractor and to
express it in terms of the physical parameters
of the system.
In particular, it may even be reduced
to a single point, thus failing in capturing the very rich and most interesting  transient behavior of
the system considered. The simplest example of such a system is the following 1D
real Ginzburg-Landau equation
$$
\Dt u=\eb\partial^2_x u+u-u^3,\ \ x\in[0,1],\ \ u\big|_{x=0}=u\big|_{x=1}=-1.
$$
In that case, the global attractor $\Cal A=\{-1\}$ is trivial for all  $\eb>0$.
However, this attractor is, factually, invisible and unreachable if
$\eb$ is small enough since the transient structures (which are very far from the attractor) have
an extremely large lifetime $T\sim e^{1/\sqrt{\eb}}$.
\par
In order to overcome these drawbacks, the concept of {\sl exponential attractor}\/ has
then been proposed in \cite{EFNT}) to possibly overcome
this difficulty. We recall below the definition of an exponential attractor adopted for
our case, see e.g., \cite{EFNT} and \cite{MZsur} for more detailed exposition.
\begin{defin}\label{Def2.ea}
A compact
subset $\mathcal{M}$ of the phase space $\Phi$ is called
an {\sl exponential attractor}\/ for the semigroup $S(t)$
if the following conditions are satisfied:\\[2mm]
{\sl (E1)}~~The set $\mathcal{M}$ is {\it positively} invariant, i.e.,
 $S(t)\mathcal{M}\subset \mathcal{M}$ for all $t\ge 0$;\\[1mm]
{\sl (E2)}~~The fractal dimension (see, e.g., \cite{mane,Te})
 of $\mathcal{M}$ in $\Phi$ is finite;\\[1mm]
{\sl (E3)}~~The set $\mathcal{M}$ attracts exponentially
 fast the image the phase space
 $\Phi$. Namely, there exist $C,\beta>0$
 such that
 \begin{equation}\label{expo-attractio}
   \dist_{L^\infty}(S(t)\Phi,\mathcal{M})\le C\+e^{-\beta t},
   \quad\perogni t\ge 0.\\[1mm]
 \end{equation}
\end{defin}
Thanks to the control of the convergence rate $(E3)$ it follows that, compared to the global attractor,
an exponential attractor
is much more robust to perturbation (usually it is H\"older continuous with respect
to the perturbation parameter, see  Section \ref{s4} below). However, since
the the exponential attractor $\Cal M$ is only {\it positively} invariant (see (E1)),
it is obviously not unique. Thus, the concrete choice of an exponential attractor
and its explicit construction becomes essential. We recall also that, in the original
paper \cite{EFNT} the construction was extremely implicit (involving the Zorn lemma) and this fact
did not allow to develop a reasonable perturbation theory. This drawback has been
overcome later in \cite{EMZ} and \cite{EMZ1} where an alternative and relatively simple
and explicit construction for the exponentially attractor has been suggested. Note also
that the construction of \cite{EMZ1}  gives an exponential attractor which
is automatically H\"older continuous with respect to the reasonable perturbations of the semigroup considered and this somehow resolves the non-uniqueness problem. We refer the reader
to the recent survey \cite{MZsur} for the detailed informations on the exponential attractors theory.
\par
The proof of the existence of an exponential attractor for our obstacle problem \eqref{eq}
will be organized as follows. First of all, we establish (in a quite standard way)
the existence of an exponential
attractor under a crucial  \emph{additional} (for this moment only) assumption on the differences of the Lagrange multiplies
$h_{u_1}$ and $h_{u_2}$ of two different solutions $u_1(t)$ and $u_2(t)$ of \eqref{eq}.
This assumption is the core of the exponential attractor
existence Theorem and it is the main result of the paper. The next
section \ref{s3} will be dedicated to its proof in two different situations:
the convex set $\K$ is smooth, or the convex set is an $n$-dimensional simplex.
\vskip5pt
\noindent
{\bf Assumption $\mathscr{L}$:} {\it There exist a closed positively invariant absorbing set $B_0\subset\Phi$
for the semigroup $S(t)$ such that, for every two solutions $u_1$ and $u_2$ of \eqref{1.eq1} starting
from $B_0$ (i.e., $u_i(0)\in B_0$), the corresponding Lagrange multipliers $h_{u_1}$
and $h_{u_2}$ satisfy the following estimates:
\begin{equation}\label{2.lag}
\|h_{u_1}-h_{u_2}\|_{L^1([0,1]\times\Omega)}\le \mathscr{C}\|u_1(0)-u_2(0)\|_{L^2},
\end{equation}
 where the constant $\mathscr{C}$ is independent of $u_1$ and $u_2$}.
\vskip5pt
Under Assumption $\mathscr{L}$, in the next Theorem we prove the existence of an exponential attractor (and thus
the finite dimensionality of the global attractor).
\begin{teorema}\label{Th2.exp} Let the assumption of Theorem \ref{Th1.wp} hold and
let, in addition, the semigroup $S(t)$ associated with equation \eqref{eq} satisfy
Assumption $\mathscr{L}$. Then, $S(t)$ possesses an exponential attractor $\Cal M\subset\Phi$ in the
sense of Definition \ref{Def2.ea}. Moreover, $\Cal M$ is a bounded subset of $C^{2-\nu}$,
for all $\nu>0$. Finally, the global attractor $\Cal A$ constructed in the Theorem \ref{Th2.globattr}
has finite fractal dimension.
\end{teorema}
\begin{proof} Recall that, since $B_0$ is a semi-invariant absorbing set for $S(t)$,
it is sufficient to verify the existence of the exponential attractor for the
restriction of $S(t)$ on $B_0$ only. As usual, we first verify the existence
of such attractor for the discrete semigroup generated by the map $S=S(1)$ and then
extend to the continuous time. To this end, we will use the following abstract
exponential attractor existence theorem suggested in \cite{EMZ}.
\begin{lemma}\label{Lem2.abs} Let $\mathscr{H}$ and $\mathscr{H}_1$ be two Banach spaces such that
$\mathscr{H}_1$ is compactly embedded to $H$. Let $\Bbb B_0$ be a bounded closed subset of $\mathscr{H}$ and
a map
$$
\Bbb S: \Bbb B_0\to\Bbb B_0
$$
be such that
\begin{equation}\label{2.sqw}
\|\Bbb S b_1-\Bbb S b_2\|_{\mathscr{H}_1}\le K\|b_1-b_2\|_{\mathscr{H}},\ \ b_1,b_2\in\Bbb B_0,
\end{equation}
where the constant $K$ is independent of $b_1$ and $b_2$. Then, the discrete
semigroup $\{\Bbb S(n),\, n\in\Bbb N\}$ generated on $\Bbb B_0$ by the iterations of the map $\Bbb S$ possesses an
exponential attractor, i.e., there exists a compact set $\Cal M_d\subset\Bbb B_0$ such that
\par
(E1) $\Cal M_d$ is positively invariant: $\Bbb S\Cal M_d\subset\Cal M_d$;
\par
(E1) The fractal dimension of $\Cal M_d$ in $\mathscr{H}$ is finite:
$$
\dim_f(\Cal M_d,\mathscr{H})\le M<\infty;
$$
and
\par
(E3) $\Cal M_d$ attracts exponentially the images of $\Bbb B_0$ under the iterations of the map $\Bbb B_0$:
$$
\dist_{\mathscr{H}}(\Bbb S(n)\Bbb B_0,\Cal M_d)\le Ce^{-kn}.
$$
Moreover, the positive constants $M$, $C$ and $k$ can be expressed explicitly in terms
of the squeezing constant $K$, the size of the set $\Bbb B_0$ and the entropy of the
compact embedding $\mathscr{H}_1\subset \mathscr{H}$.
\end{lemma}
We will use the so-called method of "$l$-trajectories (introduced by M\'alek and Ne{\v{c}}as
in \cite{MN}, see also \cite{MP} and \cite{Z}) in order to construct the proper spaces $\mathscr{H}$ and $\mathscr{H}_1$
and to verify the assumptions of Lemma \ref{Lem2.abs}.
\par
Namely, let us consider the trajectory space $\Bbb B_0$ consisting of the pieces of trajectories of the solution semigroup $S(t)$ of length one starting from $B_0$:
\begin{equation}\label{2.trphase}
\Bbb B_0:=\{u\in L^\infty([0,1],\Phi),\ u(0)=u_0\in B_0,\ \ u(t)=S(t)u_0,\ t\in[0,1]\}.
\end{equation}
Then, there is a one-to-one correspondence between $B_0$ and $\Bbb B_0$ generated
by the solution map
$$
\Bbb T\,:\,B_0\to\Bbb B_0,\ \ (\Bbb Tu_0)(t):=S(t)u_0
$$
and, therefore, we may lift the semigroup $S(t):B_0\to B_0$ to the
conjugated semigroup $\Bbb S(t)$ acting on the trajectory space $\Bbb B_0$:
\begin{equation}\label{2.lift}
\Bbb S(t):\Bbb B_0\to\Bbb B_0,\ \ \ \Bbb S(t):=\Bbb T\circ S(t)\circ\Bbb T^{-1}.
\end{equation}
We intend to apply Lemma \ref{Lem2.abs} to the map $\Bbb S=\Bbb S(1)$ acting on the
trajectory space $\Bbb B_0$. To this end, we define the spaces $\mathscr{H}$ and
 $\mathscr{H}_1$ as follows:
\begin{equation}\label{2.spaces}
\mathscr{H}:=L^2(0,1; L^2),\ \ \ \mathscr{H}_1:=L^2(0,1;H^1_{0})\cap W^{1,1}(0,1;H^{-s}),
\end{equation}
where $s>\max\{1,N/2\}$ is a fixed exponent. Then, obviously, the embedding $\mathscr{H}_1\subset \mathscr{H}$
is compact and we only need to check the smoothing property \eqref{2.sqw}. To this end,  we need (together with Assumption $\mathscr{L}$ and estimate \eqref{2.lip}) the following additional regularization property of the semigroup $S(t)$.
\begin{lemma}\label{Lem2.addreg} Let $u_1(t)$ and $u_2(t)$ be two solutions of problem
\eqref{eq}. Then, the following estimate holds:
\begin{equation}\label{2.tr-reg}
\|u_1(1)-u_2(1)\|_{L^2}^2\le (2\lambda+1)\int_0^1\|u_1(t)-u_2(t)\|_{L^2}^2\,dt.
\end{equation}
\end{lemma}
Indeed, multiplying the differential inequality \eqref{1.gr} by $t$ and
integrating $t\in[0,1]$, we arrive at \eqref{2.tr-reg}.
\par
We are now ready to verify the smoothing property \eqref{2.sqw}.
\begin{lemma}\label{Lem2.sqz} Let the assumptions of Theorem \ref{Th1.wp} hold and,
in addition, Assumption $\mathscr{L}$ be satisfied. Then the following estimate holds for every two solutions $u_1(t)$ and $u_2(t)$ such that $u_0^{i}\in B_0$, $i=1,2$:
\begin{equation}\label{2.sqz1}
\|u_1-u_2\|_{L^2(1,2;H^1_0)}+\|\Dt u_1-\Dt u_2\|_{L^1(1,2;H^{-s})}\le
L\|u_1-u_2\|_{L^2(0,1;L^2)},
\end{equation}
where the constant $L$ is independent of $u_1$ and $u_2$.
\end{lemma}
\begin{proof}
First of all we recall that, for $v\in L^1(1,2;H^{-s})$,
$$ \| v\|_{L^1(1,2;H^{-s})} = \sup_{\varphi}\Big\vert\int_{1}^{2}\;\langle v,\varphi\rangle dr\Big\vert,$$
where the $\sup$ is taken over the $\varphi \in L^\infty(1,2;H^{s}_{0})$ such that
$\|\vf\|_{L^\infty(1,2;H^{s}_{0})} = 1$ and the duality pairing is of
course between $H^{-s}$
and $H^{s}_{0}$. Consequently, there holds
\begin{eqnarray*}
&\disp \int_{1}^{2}\|\Dt u_1(t)-\Dt u_2(t)\|_{H^{-s}}dt\nonumber\\
&\disp  \le \int_{1}^{2}\|u_1(t)-u_2(t)\|_{H^1_{0}}dt+
\int_{1}^{2}\|h_{u_1}(t)-h_{u_2}(t)\|_{L^1}dt+\lambda\int_{1}^{2}\|u_1(t)-u_2(t)\|_{L^2}dt
\end{eqnarray*}
(here we have implicitly used that $s>\max\{1,N/2\}$ which implies that
$H^{s}_{0}\subset L^\infty$). Using now Assumption $\mathscr{L}$, together
with the global Lipschitz continuity \eqref{2.lip}, we have
$$
\|u_1-u_2\|_{L^2(1,2;H^1_0)}+\|\Dt u_1-\Dt u_2\|_{L^1(1,2;H^{-s})}\le
C\|u_1(1)-u_2(1)\|_{L^2}.
$$
Combining this estimate with \eqref{2.tr-reg}, we arrive at \eqref{2.sqz1} and finish
the proof of the lemma.
\end{proof}
We are now ready to finish the proof of Theorem \ref{Th2.exp}.
Indeed, we have verified that
the map $\Bbb S:=\Bbb S(1)$ satisfies all of the assumptions of the
abstract Lemma \ref{Lem2.abs} and, therefore, the semigroup
 $\{\Bbb S(n),\ n\in\Bbb N\}$ possesses
an exponential attractor $\Bbb M_d$ in the trajectory space $\Bbb B_0$ endowed with
the topology of $\mathscr{H}=L^2(0,1;L^2)$. Projecting it back to the phase space
$B_0\subset\Phi$ via
\begin{equation}\label{2.prexp}
\Cal M_d:=\Bbb M_d\big|_{t=1}
\end{equation}
and using estimate \eqref{2.tr-reg}, we see that $\Cal M_d$ is indeed the exponential
attractor for the discrete semigroup $\{S(n),\,n\in\Bbb Z\}$ acting on $B_0$ (endowed with
the topology of $L^2$). In addition, we see that
\begin{equation}\label{2.shift}
\Cal M_d\subset S(1)B_0\subset S(1)\Phi
\end{equation}
which implies that $\|\Dt u(t)\|_{L^2}\le C$ for every trajectory $u(t)$ starting from
$\Cal M_d$ (due to estimate \eqref{1.reg}). Thus, thanks to \eqref{2.lip}, the map
$(t,u_0)\to S(t)u_0$ is globally Lipschitz continuous on $[0,1]\times\Cal M_d$ (in the
$\R\times L^2(\Omega)$-metric). Thus, the
desired exponential attractor for the solution semigroup $S(t)$ with {\it continuous} time
can be now constructed by the following standard expression (see \cite{EFNT} for
the details):
$$
\Cal M:=\cup_{t\in[0,1]}\Cal M_d.
$$
Note that, up to now, we have verified that $\Cal M$ has the finite fractal dimension
and possesses the exponential attraction property in the topology of $L^2$ only.
In order to verify that these properties actually hold in the $L^\infty$-topology
of the phase space $\Phi$, we use (in a standard way) the additional regularity of
$\Cal M$ and a proper interpolation inequality. Indeed, by our construction,
$\Cal M\subset S(1)\Phi$ and, therefore, due to estimate \eqref{1.c2}, $\Cal M$ is
globally bounded in $C^{2-\nu}$. Using now the following interpolation inequality
\begin{equation}\label{2.int}
\|u-v\|_{L^\infty}\le C_\nu\|u-v\|^\kappa_{L^\infty}\|u-v\|^{1-\kappa}_{C^{2-\nu}},\ 0<\kappa<1,
\end{equation}
we see that the dimension of $\Cal M$ is finite not only in $L^2$, but also in $L^\infty$ and that
the attraction property holds in $L^\infty$ as well.
Thus, the desired exponential attractor
$\Cal M$ in the phase space $\Phi$ is constructed and Theorem \ref{Th2.exp} is proved.
\end{proof}
\begin{remark}\label{rem2.standard} Note once more that the method of proving the
{\it conditional} result of Theorem \ref{Th2.exp} (which is  more or less standard variation
of the $l$-trajectories method)
is widely used now-a-days in the exponential attractors theory, see e.g., \cite{MZsur}
and the references therein. Thus, the major difficulty here (and the major novelty of the paper)
is related with the verification of the Assumption $\mathscr{L}$ for the solution semigroup $S(t)$
associated with the obstacle problem \eqref{eq}.
\end{remark}

\section{Estimates on the difference of the Lagrange multipliers and verification of
Assumption $\mathscr{L}$}\label{s3}
The aim of this section is to verify the Assumption $\mathscr{L}$ under some additional assumptions on the
structure of the convex set $\mathscr{K}$. We start with the case in which
$\mathscr{K}$ has a smooth boundary.

\subsection{The case of regular $\mathscr{K}$.}  The main result of this subsection is the following Theorem.
\begin{teorema}\label{Th3.smooth} Let the assumptions of Theorem \ref{Th1.wp} hold and let,
in addition, the boundary $S=\partial \mathscr{K}$ be smooth enough (at least, $C^{2,1}$). Then, the solution semigroup
$S(t)$ associated with equation \eqref{eq} satisfies  Assumption $\mathscr{L}$.
\end{teorema}
\begin{proof} The proof is based on an argument that combines a proper choice of an approximation scheme
and the maximum principle  similar to the one devised to prove Proposition \ref{Prop1.linf}.
However, since the function $M(u)=\hbox{dist}(u,\mathscr{K})$
used to define $F_\epsi$ in
the existence Theorem \ref{Th1.wp} is not smooth enough near the boundary,
we should introduce another approximation of the singular potential.
The smooth
correction of $M$ is given by the following Lemma.
\begin{lemma}\label{Lem3.msmooth} Let $\mathscr{K}$ be a bounded convex set with the $C^{2,1}$-smooth boundary $S$.
Then, there exists
a function $M:\mathbb{R}^n\to \mathbb{R} $
with the following properties:
\begin{subequations}\label{1}
\begin{align}
&\disp  M\in C^{2,1}(\R^n);\label{m1}\\
&\disp  M\ \text{is convex};\label{m2}\\
&\disp S=\partial \mathscr{K}=\left\{z\in \Bbb R^n \hbox{ such that }M(z)=0\right\};\label{m3}\\
&\disp |\nabla M(z)|=\theta(M(z))\label{m4},
\end{align}
\end{subequations}
where $\teta=\theta(z)$ is a monotone increasing function which is smooth near $z=0$
and such that  $\theta(0)\ne0$.
\end{lemma}
\begin{proof}
The function $M$ can be constructed as follows.
First of all, for any $\delta>0$, we introduce the following set:
\begin{equation}\label{Sdelta}
S_{-\delta}:=\left\{z_0\in \cK: \hbox{ dist}(z_0,S)= \delta\right\}.
\end{equation}
Then, we consider the subset $\mathscr{K}_{-\delta}\subset \mathscr{K}$
whose boundary is $S_{-\delta}$.
It turns out that, being $\mathscr{K}$ convex, the domain $\mathscr{K}_{-\delta}$ is convex too. Moreover,
by possibly taking $\delta$ small, the boundary $S_{-\delta}$ has the same regularity of $S$.
The candidate for $M$ is thus the function
\begin{equation}\label{2}
M(z):=[\dist(z,\mathscr{K}_{-\delta})]^3-\delta^{3}.
\end{equation}
In fact, $M$ is clearly convex and regular. Moreover, by choosing $\teta=\teta(w)=3 (w+\delta^3)^{2/3}$ and using the identity
$$
|\nabla\left(\dist(z,\mathscr{K}_{-\delta})\right)|=1,
$$
one sees that the condition \eqref{m4} is satisfied.

Finally, by possibly taking a small
$\delta$, the smoothness of the boundary $S$ entails the validity of \eqref{m3}.
Note also that $M$ has the additional property of qualifying the fact that
$v\in \mathscr{K}$, namely it turns out
that  $v\in \mathscr{K}$ if and only if $-\delta^{3}\le M(v) \le 0$.
\end{proof}
Now, we introduce the approximations $F_\eb(u)$ of the indicator
function $I_{\mathscr{K}}$ as in \eqref{1.dddd} using the above defined function $M$.
Let now $u(t)$ and $v(t)$ be two solutions of the singular
equation \eqref{1.eq1} with initial conditions
$$
u_0,v_0 \in C^{2-\nu}\cap H^1_0,
$$
with sufficiently small positive $\nu$ (actually, $\nu=1$ is sufficient for what follows).
Then, let $u_\eb(t)$ and $v_\eb(t)$ be their approximations, namely the solutions of  the approximate
problems \eqref{1.eqreg} with initial data $u_0$ and $v_0$.
As one can easily check,  Proposition \ref{Prop1.linf} and
Corollary \ref{Cor1.c2} still hold for this new approximation entailing the validity of the following two estimates
\begin{equation}\label{3.uni}
\|u_\eb(t)\|_{C^{2-\nu}}+\|v_\eb(t)\|_{C^{2-\nu}}\le Q(\|u_0\|_{C^{2-\nu}}+\|v_0\|_{C^{2-\nu}}),\ \ t\ge0,
\end{equation}
and
\begin{equation}\label{3.unibis}
\|F'_\eb(u_\eb(t))\|_{L^\infty}+\|F'_\eb(v_\eb(t))\|_{L^\infty}\le C,
\end{equation}
where the monotone function $Q$ and the constant $C$ are both independent of $\eb$. Thus, up to not relabeled subsequence,
$$
h_{u_\eb}:=F'_\eb(u_\eb)\to h_u,\ \ h_{v_\eb}:=F'_\eb(v_\eb)\to h_v
$$
weakly strar in $L^\infty([0,T]\times\Omega)$.
Thus, in order to prove the theorem, we only need
to verify that
\begin{equation}\label{3.heps}
\|F'_\eb(u_\eb)-F'_\eb(v_\eb)\|_{L^1([0,1]\times\Omega)}\le Q(\|u_0\|_{C^{2-\nu}}+\|v_0\|_{C^{2-\nu}})\|u_0-v_0\|_{L^2}
\end{equation}
with the function $Q$ not depending on $\eb$ and then by semicontinuity we obtain the desired estimate \eqref{2.lag} in the limit $\epsi\searrow 0$.
This estimate will hold for every two trajectories $u$ and $v$ such that $u_0$ and $v_0$ are bounded in $C^{2-\nu}$.
In order to construct the desired positively invariant absorbing set $B_0$, it is sufficient to
take a ball
$$
B:=\{u\in C^{2-\nu/2}\cap H^1_0\cap \Phi,\ \ \|u\|_{C^{2-\nu}}\le R\}
$$
with a sufficiently large radius $R$. Then, thanks to Corollary \ref{Cor1.c2},
 $B$ will be an {\it absorbing}
set for the semigroup $S(t)$ associated  with equation \eqref{eq} and, therefore,
the closure $B_1:=[B]_{\Phi}$ of the set $B$ in the topology of $\Phi$ will be a bounded in $C^{2-\nu}$
 and {\it closed} in $\Phi$ absorbing set for the semigroup $S(t)$. Applying Corollary \ref{Cor1.c2} again,
we see that the set
\begin{equation}\label{3.b0}
B_0:=\bigcup_{t\ge0}B_1
\end{equation}
will be the desired positively invariant, closed in $\Phi$ and bounded in $C^{2-\nu}(\Omega)$ absorbing set
for semigroup $S(t)$ and estimate \eqref{3.heps} will guarantee that \eqref{2.lag} will hold uniformly with respect to
all trajectories starting from $B_0$.
\par
Thus, it only remains to verify the uniform estimate \eqref{3.heps}.
The first step is, roughly speaking, to reduce \eqref{eq} to an equation with scalar constraint
(actually to an approximation of). This is done using the function $M$ as we did in \eqref{1.max}. Indeed,
testing equation \eqref{1.eqreg} in the scalar product $(\cdot,\cdot)_{\mathbb{R}^n}$ of
$\mathbb{R}^n$ with $\nabla M(u)$ and using the condition \eqref{m4}, we obtain
\begin{equation}\label{7}
\Dt(M(u_\eb)) - \Delta(M(u_\eb)) + \mathcal{F}'_\eb(M(u_\eb)) + D(u_\eb) = 0,\
\end{equation}
with $\mathcal{F}'_\eb(z):=f'_\eb(z)\theta^2(z)$
and
\begin{equation}\label{8}
D(u):=\sum_{i=1}^{N}\sum_{j,k=1}^{n}M''_{j,k}(u)\partial_{x_i}u_j\partial_{x_i}u_k-\lambda\sum_{j}M'_j(u)u_j,
\end{equation}
where $\left\{M_{j,k}\right\}_{j,k=1}^{n}$ denotes the entries of the Hessian matrix $\Cal H(M)$ of the function
  $M$ (and the analogous equation holds for $v_\eb$).
Let now $w_\eb:=M(u_\eb)-M(v_\eb)$. Then
\begin{equation}\label{9}
\Dt w_\eb -\Delta w_\eb  + [\mathcal{F}'_\eb(M(u_\eb))-\mathcal{F}'_\eb(M(v_\eb))] + [D(u_\eb)-D(v_\eb)] = 0.
\end{equation}
We now multiply this equation by $\sgn w_\epsi$ and integrate over $\Omega$ obtaining
\begin{eqnarray}\label{9bis}
&\disp \frac{d}{dt} \|u_\eb-v_\eb\|_{L^{1}} + \|\mathcal{F}'_\eb(M(u_\eb))-\mathcal{F}'_\eb(M(v_\eb))\|_{L^{1}}
\nonumber\\
&\disp =  -
\int_{\Omega}(D(u_\eb)-D(v_\eb)) \sgn w_\epsi dx =:\mathcal{R},
\end{eqnarray}
where we have used the monotonicity of $\mathcal{F}'_\eb(z)$ (recall
in particular that $\teta(z)= (z+\delta^3)^{2/3}$) and the Kato
inequality. To estimate the term with $\mathcal{R}$ we start to
rewrite it in the more explicit form
\begin{eqnarray}\label{R}
&\disp \mathcal{R} =
\sum_{i=1}^{N}\int_\Omega \left((\mathcal{H}(M(u_\eb))-(\mathcal{H}(M(v_\eb)))\partial_{x_i}u_\eb,
\partial_{x_i}u_\eb\right)_{\mathbb{R}^n}\sgn(w_\eb) dx \nonumber\\
&\disp + \sum_{i=1}^{N}\int_\Omega\left(\mathcal{H}(M(v_\eb))(\partial_{x_i}u_\eb-\partial_{x_i}v_\eb),
\partial_{x_{i}}u_\eb\right)_{\mathbb{R}^n}\sgn(w_\eb)dx\nonumber\\
&\disp + \sum_{i=1}^{N}\int_\Omega\left(\mathcal{H}(M(v_\eb))(\partial_{x_i}u_\eb-\partial_{x_i}v_\eb),
\partial_{x_{i}}v_\eb\right)_{\mathbb{R}^n}\sgn(w_\eb)dx\nonumber\\
&\disp +\lambda\int_\Omega \left(\nabla M(u_\eb),u_\eb\right)_{\mathbb{R}^n}\sgn(w_\eb) dx -
\lambda\int_\Omega \left(\nabla M(v_\eb),v_\eb\right)_{\mathbb{R}^n}\sgn(w_\eb) dx.
\end{eqnarray}
 Thus, thanks to
the assumed $C^{2,1}$ regularity on the function $M$ (which follows from the
smothness of the boudary of $\mathscr{K}$) and the $C^1$-regularity of $u_\eb$ and $v_\eb$
which follows from \eqref{3.uni},
it is not difficult to realize that
\begin{equation}\label{Rbis}
\|\mathcal{R}(t)\|_{L^1(\Omega)} \le C\|u_\eb(t)-v_\eb(t)\|_{H^1} \hbox{ for any } t\ge 0,
\end{equation}
where the constant $C$ depends on the $C^1$-norms of $u_0$ and $v_0$, but is independent of $\eb>0$
\par
Integrating the differential inequality \eqref{9bis} with respect to $t\in[0,1]$ and using \eqref{Rbis}
and the analog of \eqref{2.lip} for the solutions $u_\eb$ and $v_\eb$ of the approximate problems \eqref{1.eqreg}, we arrive at
(recall also that $u_\epsi(0) = u_0$ and $v_\epsi(0) = v_0$)
\begin{eqnarray}\label{3.fest}
&\disp \int_0^1\|\Cal F'_\eb(M(u_\eb))-\Cal F'_\eb(M(v_\eb))\|_{L^1}\,dt\le \|u_0-v_0\|_{L^1} +
\int_{0}^{1}\Cal R(t)dt
\nonumber\\
&\disp   \le \|u_\eb(0)-v_\eb(0)\|_{L^1}+C\int_0^1\|u_\eb(t)-v_\eb(t)\|_{H^1}\,dt\nonumber\\
&\disp \le
C_1\|u_0-v_0\|_{L^2},
\end{eqnarray}
where the constants $C$ and $C_1$ depend on the $C^1$-norms of $u(0)$ and $v(0)$, but are independent of $\eb$.
\par
The next step is to estimate the difference $f'_\eb(M(u_\eb))-f'_\eb(M(v_\eb))$ in terms of the difference
of $\Cal F'(M(u_\eb))$ and $\Cal F'_\eb(M(v_\eb))$. To this end, we note that $f_\eb(z)\equiv0$ if $z\le0$, so, without loss
of generality, we may assume that, say, $M(u_\eb)\ge0$. Then, if $M(v_\eb)\le0$
\begin{eqnarray}\label{3.Ff}
&\disp |f'_\eb(M(u_\eb))-f'_\eb(M(v_\eb))|=f'_\eb(M(u_\eb))=\theta^{-2}(M(u_\eb))\Cal F'_\eb(M(u_\eb))\nonumber\\
&\disp \le
\theta^{-2}(0)\Cal F'_\eb(M(u_\eb))=C|\Cal F'_\eb(M(u_\eb))-\Cal F'_\eb(v_\eb)|.
\end{eqnarray}
Let now  $M(v_\eb)\ge0$. Then,
\begin{eqnarray}\label{3.fF}
&\disp |f'_\eb(M(u_\eb))-f'_\eb(M(v_\eb))|= |\Cal F'_\eb(M(u_\eb))\theta^{-2}(M(u_\eb))-\Cal F'_\eb(M(v_\eb))\theta^{-2}(M(v_\eb))|\nonumber\\
&\disp \le \theta^{-2}(M(u_\eb))\cdot|\Cal F'_\eb(M(u_\eb))-\Cal F'_\eb(M(v_\eb))|+
\Cal F'_\eb(M(v_\eb))|\theta^{-2}(M(u_\eb))-\theta^{-2}(M(v_\eb))|\nonumber \\
&\disp\le \theta^{-2}(0)|\Cal F'_\eb(M(u_\eb))-\Cal F'_\eb(M(v_\eb))|+
\|F'_\eb(u_\eb)\|_{L^\infty}\theta^{-3}(0)|\theta^2(M(u_\eb))-\theta^2(M(v_\eb))|\nonumber \\
&\disp \le C\Big(|\Cal F'_\eb(M(u_\eb))-\Cal F'_\eb(M(v_\eb))|+|u_\eb-v_\eb|\Big),
\end{eqnarray}
where the constant $C$ is independent of $\eb$ (here we have implicitly used that the $L^\infty$-norms
of $F'_\eb(u_\eb)$ are uniformly bounded, thanks to \eqref{3.unibis}). Thus, due to \eqref{3.fest},
\eqref{3.fF} and \eqref{3.Ff}
$$
\int_0^1\|f'_\eb(M(u_\eb(t)))-f'_\eb(M(v_\eb(t)))\|_{L^1}\,dt\le C\|u_\eb(0)-v_\eb(0)\|_{L^2}.
$$
Finally, since $M$ is at least $C^2$,
\begin{eqnarray}\label{11}
&\disp \|F'_\eb(u_\eb)-F'_\eb(v_\eb)\|_{L^1}=\|f'_\eb(M(u_\eb)M'(u_\eb)-f'_\eb(M(v_\eb))M'(v_\eb)\|_{L^1}\nonumber\\
&\disp \le
\|M'(u_\eb)\|_{L^\infty}\|f'_\eb(M(u_\eb))-f'_\eb(M(v_\eb)\|_{L^1}+
\|f'_\eb(M(u_\eb))\|_{L^\infty}\|M'(u_\eb)-M'(v_\eb)\|_{L^1}\nonumber\\
&\disp \le
C(\|f'_\eb(M(u_\eb))-f'_\eb(M(v_\eb))\|_{L^1}+\|u_\eb-v_\eb\|_{L^1})
\end{eqnarray}
and
$$
\int_0^1\|F'_\eb(M(u_\eb(t))-F'(M(v_\eb(t))\|_{L^1}\,dt\le C\|u_\eb(0)-v_\eb(0)\|_{L^2}.
$$
Thus, estimate \eqref{3.heps} is verified and Theorem \ref{Th3.smooth} is proved.
\end{proof}

\begin{coro}\label{Cor3.exp} Let the assumptions of Theorem \ref{Th1.wp} hold and
let boundary $S=\partial \mathscr{K}$ be of class $C^{2,1}$. Then, the semigroup $S(t)$ associated
with the obstacle problem \eqref{eq} possesses an exponential attractor $\Cal M$ in the sense
of Definition \ref{Def2.ea} in the phase space $\Phi$.
Moreover, the global attractor has finite fractal dimension.
\end{coro}
Indeed, this assertion is an immediate corollary of Theorems \ref{Th2.exp} and \ref{Th3.smooth}.

\subsection{The case of an irregular convex $\mathscr{K}$: the simplex \eqref{phasesimplex}.}
As we  saw in the former section,  the
the smoothness of the boundary $S=\partial \mathscr{K}$  seems crucial for the method of verifying
 Assumption $\mathscr{L}$ suggested in the proof of Theorem \ref{Th3.smooth}.
 First of all,  the $C^{2,1}$ smoothness of $M$ is necessary in order to obtain estimate \eqref{Rbis}. Secondly, the fact that $\theta(0)\ne 0$ in \eqref{m4}, which is crucial to obtain estimates \eqref{3.fF} and \eqref{3.Ff}) implies that $\nabla M(z)\ne0$ for all $z\in S$. Thus, by the implicit function theorem, the boundary $S$ also must be at least $C^{2,1}$-smooth.
However, from the possible applications to phase transition problems, one of the most important examples for the set $\mathscr{K}$ is the
$n$-dimensional simplex, namely the set
\begin{equation}\label{3.phasesimplex}
\mathscr{K}:= \left\{p=(p_1,\ldots,p_n)\in \mathbb{R}^n \mbox{ such that } \sum_{i=1}^{n}p_1\le 1,\ \ p_i\ge0,\ i=1,\cdots,n\right\}
\end{equation}
which is a polyhedron  and its boundary is only {\it piece-wise} smooth. Thus, Theorem \ref{Th3.smooth} is not directly
applicable here. Nevertheless, as we will show
below, Assumption $\mathscr{L}$ remains valid for  the non-regular case \eqref{3.phasesimplex}. Again, the verification of Assumption $\mathscr{L}$ will rely
on an approximation argument and, to this purpose, we will consider a slightly different
approximation for the indicator function $I_{\mathscr{K}}$, namely, let us introduce
\begin{equation}\label{3.appr}
F_\eb(u):=f_\eb(\sum_{i=1}^n u_i-1)+\sum_{i=1}^nf_\eb(-u_i),
\end{equation}
where the function $f_\eb$ is defined by \eqref{1.ddd}. Obviously, $F_\eb(u)$ is convex
and satisfies the conditions of Remark \ref{Rem1.app}. Therefore, we can indeed use it for approximating the singular problem \eqref{eq}. Consequently, all of the estimates and convergences obtained in the proof of Theorem \ref{Th1.wp} hold for this new approximation.
\par
On the other hand, an inspection in the proof of the $L^\infty$ estimate in Proposition \ref{Prop1.linf} reveals that one of keys
point was the structure of the approximations. It is also evident that the new defined approximation for the simplex is structurally different from the approximation
introduced before; thus the $L^\infty$ bound of Proposition \ref{Prop1.linf}  need a different proof.  In particular, having in mind also the study of the approximation
of the exponential attractor (see the next section), we formulate the analogue of Proposition \ref{Prop1.linf} in a slightly stronger form by indicating also the
invariant regions at the $\epsi$-level.
Thus, the global bound result on the Lagrange multipliers takes this form
\begin{teorema}\label{Th3.linf} Let $\eb>0$ be small enough and let the function $F_\eb(u)$ be defined by \eqref{3.appr}. Then, there exists a positive constant $p$ (independent of $\eb$) such that the set
\begin{equation}\label{3.inv}
\mathscr{K}_\eb:=\{u\in\R^n,\ \ F_\eb(u)\le p\eb\}
\end{equation}
 is an invariant region for the solution semigroup $S_\eb(t)$
of the approximate problems \eqref{1.eqreg}, namely
\begin{equation}\label{3.inv1}
S_\eb(t):\Phi_\eb\to\Phi_\eb,\ \ t\ge0,
\end{equation}
where
\begin{equation}\label{3.phase}
\Phi_\eb:=\{u\in L^\infty(\Omega),\ \ u(x)\in \mathscr{K}_\eb\ \ \text{ for almost all } x\in\Omega\}.
\end{equation}
Finally, the approximations $h_{u_\eb}(t):=F'_\eb(u_\eb(t))$ (where $u_\eb(t):=S_\eb(t)u_0$, $u_0\in\Phi_\eb$) are uniformly bounded in the $L^\infty$-norm:
\begin{equation}\label{3.linf}
\|h_{u_\eb}(t)\|_{L^\infty(\Omega)}\le C
\end{equation}
where the constant $C$ is independent of $\eb$ and of the concrete choice of $u_0\in\Phi_\eb$.
\end{teorema}
\begin{proof}First of all, we need the following Lemma which clarifies the relations
between $|F'_\eb(u)|$ and $F_\eb(u)$.
\begin{lemma}\label{Lem3.F} Let the function $F_\eb$ be defined by \eqref{3.appr} and
the function $f_\eb(z)$ be given by \eqref{1.ddd}. Then, there exist two positive constants $\kappa_1$ and $\kappa_2$ (independent of $\eb$) such that
\begin{equation}\label{3.fpr}
\frac{\kappa_2}\eb F_\epsi(u)\le |F'_\eb(u)|^2\le\frac{\kappa_1}\eb F_\eb(u)
\end{equation}
for all $\eb>0$ and all $u\in\R^n$.
\end{lemma}
\begin{proof} Indeed, let $f'_i:=-f'_\eb(-u_i)$ and $f'_0:=f'_\eb(\sum u_i-1)$. Then,
\begin{equation}\label{3.est}
|F'_\eb(u)|^2=\sum_{i=1}^n(f'_{i}+f'_{0})^2=n \vert f'_{0}\vert ^2+\sum_{i=1}^n \vert f'_{i}\vert^2+2(\sum_{i=1}^n f'_{i})f'_{0}.
\end{equation}
Let us consider two cases:
\par
\noindent
Case I: $f'_0=f'_{\eb}(\sum u_i-1)=0$. Then, \eqref{3.est} simply reads as
$$
|F'_\eb(u)|^2=\sum_{i=1}^n \vert f'_{i}\vert ^2=\sum_{i=0}^{n}\vert f'_{i}\vert^2.
$$
\par
\noindent
Case II: $f'_0=f'_\eb(\sum u_i-1)\ne0$. Then, keeping in mind the definition of
the function $f_\eb$, we conclude that at least one of $u_i$, $i=1,\cdots,n$ (say,
$u_n$ for definiteness) must be positive and, therefore, without loss of generality,
we  may assume that $u_n>0$ and, consequently, $f'_n=-f'_\epsi(-u_n)=0$. Then, using the elementary
inequality
$ab\ge-\frac12(\frac {a^2}\alpha+\alpha b^2)$, $\alpha>0$, we arrive at
$$
n\vert f'_0\vert^2+\sum_{i=1}^{n-1}\vert f_i\vert^2+\sum_{i=1}^{n-1}f'_0f'_i\ge (n-\alpha(n-1))\vert f'_0\vert^2+(1-\alpha^{-1})\sum_{i=1}^{n-1} \vert f'_i\vert^2.
$$
Choosing the positive $\alpha=\alpha(n)$ in an optimal way as a solution of the equation
$$
1-\alpha^{-1}=n(1-\alpha)+\alpha),
$$
we see that, in the second case
$$
|F'_\eb(u)|^2\ge\theta(n)\sum_{i=0}^n \vert f'_i\vert^2
$$
with $\theta(n):=\frac{n+1-\sqrt{(n+1)^2-4}}4>0$. Thus, in both cases
$$
|F'_\eb(u)|^2\ge\theta(n)\sum_{i=0}^n \vert f'_i\vert^2 = \theta(n)\Big(\vert f'_0\vert^2 + \sum_{i=1}^n\vert f'_i\vert^2\Big).
$$
Since the upper bound is obvious, we arrive at the following inequality
\begin{equation}\label{3.est1}
\theta(n)\(\sum_{i=1}^n |f'_\eb(-u_i)|^2+|f'_\eb(\sum_{i=1}^n u_i-1)|^2\)\le |F'_\eb(u)|^2\le
\sum_{i=1}^n |f'_\eb(-u_i)|^2+|f'_\eb(\sum_{i=1}^n u_i-1)|^2.
\end{equation}
It only remains to note that, due to our choice \eqref{1.ddd} of the function $f_\eb$,
$$
|f'_\eb(z)|^2=\frac 4\eb f_\eb(z)
$$
and, consequently, \eqref{3.est1} implies \eqref{3.fpr} and finishes the proof of the Lemma.
\end{proof}
With the help of Lemma \ref{Lem3.F}, it is now not difficult to complete the proof of Theorem \ref{Th3.linf}. Indeed, testing
equation \eqref{1.eqreg} in the scalar product of $\mathbb{R}^n$ with
$F'_\eb(u)$ and using that $F''_\eb(u)\ge0$, we have
(compare with \eqref{1.max})
$$
\Dt F_\eb(u)-\Dx(F_\eb(u))+|F'_\eb(u)|^2\le\lambda (u,F'_\eb(u))_{\mathbb{R}_n}.
$$
Using now Lemma \ref{Lem3.F} together with the obvious fact that
$$
|u|\le F_\eb(u)+C,
$$
where $C$ is independent of $\eb\to0$, we arrive at the differential inequality
for the scalar function $V(t):=F_\eb(u(t))$:
\begin{equation}\label{3.sc}
\Dt V-\Dx V+(\frac{\kappa_2}{2\eb}-\lambda^2) V\le C_1, \ \ V\big|_{\partial\Omega}=0,
\end{equation}
where the positive constants $C_1$ and $\kappa$ are independent of $\eb$ and the concrete
choice of the solution $u$.
Applying the comparison principle for the heat equations to \eqref{3.sc}, we see that
the region $\{u,\  V(u)\le p_\eb\}$ will be invariant if
$$
p_\eb\ge\frac{C_1}{\frac{\kappa_2}{2\eb}-\lambda^2}=\frac{2C_1\eb}{\kappa_2-\lambda^2\eb}.
$$
Thus, taking $p:=\frac{3C_1}{\kappa_2}$, we see that the region $\mathscr{K}_\eb$ defined
by \eqref{3.inv} will be indeed invariant with respect to the semi-flow $S_\eb(t)$
if $\eb>0$ is small enough. Thus, we only need to check estimate \eqref{3.linf}.
To this end, it remains to observe that, due to Lemma \ref{Lem3.F}
$$
|F'_\eb(u(t))|^2\le\frac{\kappa_1}\eb F_\eb(u(t))\le \frac{\kappa_1}{\eb}\cdot p\eb=\kappa_1p.
$$
Theorem \ref{Th3.linf} is proved.
\end{proof}
We are now ready to formulate the main result of this subsection.
\begin{teorema}\label{Th3.simplex} Let the assumptions of Theorem \ref{Th3.linf} hold and $\eb>0$
be small enough. Then, for any two solutions $u_\eb(t):=S_\eb(t)u_0$ and
$v_\eb(t):=S_\eb(t)v_0$ of the approximate problems \eqref{1.eqreg} starting from
$\Phi_\eb$ ($u_0,v_0\in\Phi_\eb$), the associated Lagrange multipliers
$h_{u_\eb}(t):=F'_\eb(u_\eb(t))$ and $h_{v_\eb}(t):=F'_\eb(v_\eb(t))$ satisfy the following estimate:
\begin{equation}\label{3.Lest}
\int_0^1\|h_{u_\eb}(t)-h_{v_\eb}(t)\|_{L^1(\Omega)}\,dt\le L\|u_0-v_0\|_{L^2(\Omega)},
\end{equation}
where the constant $L$ is independent of $\eb$, $u_0$ and $v_0$. In particular,
the limit senigroup $S(t)$ associated with the obstacle problem \eqref{eq} satisfies
Assumption $\mathscr{L}$ with $B_0:=\Phi$.
\end{teorema}
\begin{proof}
In order to simplify the notations
we forget for a while the $\epsi$ dependence and will write $u(t):=S_\eb(t)u_0$
and $v(t):=S_\eb(t)v_0$. Then, the function $w(t):=u(t)-v(t)$ solves
\begin{equation}\label{3.eqreg}
\Dt w_i-\Dx w_i+[f'_\eb(-v_i)-f'_\eb(-u_i)]+[f'_\eb(\sum_{i=1}^n u_i-1)-f'_\eb(\sum_{i=1}^n v_i-1)]=
\lambda w_i,\ \ i=1,\cdots, n.
\end{equation}
Multiply now the $i$th equation of \eqref{3.eqreg} by $\sgn w_i$ and integrate over
$\Omega$. Then, taking the sum over $i$ and using the Kato inequality, we arrive at
\begin{equation}
\Dt \|w\|_{L^1}-\lambda\|w\|_{L^1}+\sum_{i=1}^n\|f_\eb'(-u_i)-f'_\eb(-v_i)\|_{L^1}\le
-\big([f'_\eb(\sum_{i=1}^n u_i-1)-f'_\eb(\sum_{i=1}^n v_i-1)],\sum_{i=1}^n \sgn w_i\big).
\end{equation}
Let us estimate the right-hand side of this inequality. To this end, we first note that
$$
|\sum_{i=1}^n\sgn w_i|\le n.
$$
In addition, in the worst case where the equality holds, all $w_i$ are of the same sign.
Therefore, at that point the term $f_\eb'(\sum_{i=1}^n u_i-1)-f'_\eb(\sum_{i=1}^n v_i-1)$ will also have
the same sign and its product with $\sum_{i=1}^n\sgn w_i$ will be non-negative. Thus,
$$
-[f'_\eb(\sum_{i=1}^n u_i-1)-f'_\eb(\sum_{i=1}^n v_i-1)]\cdot\sum_{i=1}^n\sgn w_i\le (n-1)|f'_\eb(\sum_{i=1}^n u_i-1)-f'_\eb(\sum_{i=1}^n v_i-1)|
$$
and
\begin{equation}\label{3.eq1}
\Dt\|w\|_{L^1}+\sum_{i=1}^n\|f'_\eb(-u_i)-f'_\eb(-v_i)\|_{L^1}\le \lambda\|w\|_{L^1}+(n-1)\|f'_\eb(\sum_{i=1}^n u_i-1)-f'_\eb(\sum_{i=1}^n v_i-1)\|_{L^1}.
\end{equation}
Next, in order to estimate the right-hand side of \eqref{3.eq1}, we sum all the equations
\eqref{3.eqreg} and multiply the obtained relation by
$\sgn(\sum_{i=1}^n u_i-\sum_{i=1}^n v_i)$. Then, using again the Kato inequality, we have
\begin{eqnarray}\label{3.eq2}
&\disp \Dt\|\sum_{i=1}^n w_i\|_{L^1}+n\|f'_\eb(\sum_{i=1}^n u_i-1)-f'_\eb(\sum_{i=1}^n v_i-1)\|_{L^1}\nonumber\\
&\disp \le\lambda\|\sum_{i=1}^n w_i\|_{L^1}+\sum_{i=1}^n\|f'_\eb(-u_i)-f'_\eb(-v_i)\|_{L^1}.
\end{eqnarray}
Multiplying inequality \eqref{3.eq2} by $\frac n{n+1}$ and summing it to inequality
\eqref{3.eq1}, we finally arrive at
\begin{eqnarray}\label{3.eq2bis}
&\disp \Dt\(\|w\|_{L^1}+\frac n{n+1}\|\sum_{i=1}^n w_i\|_{L^1}\)\nonumber\\
&\disp +\frac1{n+1}\(\|f'_\eb(\sum_{i=1}^n u_i-1)-f'_\eb(\sum_{i=1}^n v_i-1)\|_{L^1}+
\sum_{i=1}^n\|f'_\eb(-u_i)-f'_\eb(-v_i)\|_{L^1}\)\nonumber\\
&\disp \le \lambda\(\|w\|_{L^1}+\frac n{n+1}\|\sum_{i=1}^n w_i\|_{L^1}\).
\end{eqnarray}
Applying the Gronwall inequality to this relation, we infer
\begin{eqnarray}\label{3.est3}
&\disp \int_0^1\|f'_\eb(\sum_{i=1}^n u_i(t)-1)-f'_\eb(\sum_{i=1}^n v_i(t)-1)\|_{L^1}+
\sum_{i=1}^n\|f'_\eb(-u_i(t))-f'_\eb(-v_i(t))\|_{L^1}\,dt\nonumber\\
&\disp \le
(n+1)e^\lambda\(\|w(0)\|_{L^1}+\frac n{n+1}\|\sum_{i=1}^n w_i(0)\|_{L^1}\)\le C\|u(0)-v(0)\|_{L^2}
\end{eqnarray}
which together with the obvious inequality
$$
\|h_{u_\eb}(t)-h_{v_\eb}(t)\|_{L^1}\le C\(\|f'_\eb(\sum_{i=1}^n u_i-1)-f'_\eb(\sum_{i=1}^n v_i-1)\|_{L^1}+
\sum_{i=1}^n\|f'_\eb(-u_i)-f'_\eb(-v_i)\|_{L^1}\)
$$
finishes the proof of estimate \eqref{3.Lest}. Passing now to the limit $\eb\to0$
in that estimate (analogously to the proof of Theorem \ref{Th3.smooth}), we see that
the limit semigroup $S(t)$ generated by the obstacle problem \eqref{eq} satisfies indeed
 Assumption $\mathscr{L}$ with $B_0:=\Phi$. Theorem \eqref{Th3.simplex} is proved.
\end{proof}

\begin{coro}\label{Cor3.exp1} Let the assumptions of Theorem \ref{Th1.wp} hold and let
the set $\mathscr{K}$ be an $n$-dimensional simplex \eqref{3.phasesimplex}. Then, the semigroup
$S(t)$ associated with the obstacle problem \eqref{eq} possesses an exponential attractor
$\Cal M$ in the sense of Definition \ref{Def2.ea}. Moreover, the global attractor $\Cal A$ has
finite fractal dimension.
\end{coro}
Indeed, this assertion is an immediate corollary of Theorems \ref{Th2.exp} and \ref{Th3.simplex}.

\subsection{Some generalizations} We now discuss the applications of our method to more general problems.
We start with the obvious observation that all the above results remain valid if we replace the term $\lambda u$ in the
left-hand side of equation \eqref{eq} by any  sufficiently regular interaction
function $g(u,x)$. Namely, consider
the problem
\begin{eqnarray}\label{3.eq}
\partial_t u -\Delta_{x} u + \partial I_{K}(u) +g(x,u)\ni 0,
\end{eqnarray}
where $g\in C(\Omega,C^1(\R^n,\R^n))$ is an arbitrary interaction function. Then, the following result holds.

\begin{teorema}\label{Th3.intr} Let $\mathscr{K}$ be a convex bounded set of $\R^n$ containing zero with a
 smooth boundary
(or let $K$ be an $n$-dimensional simplex \eqref{3.phasesimplex}) and let $g\in C(\Omega,C^1(\R^n,\R^n))$ be
an arbitrary (not necessarily a gradient!) non-linear interaction function. Then, the solution semigroup $S(t)$
associated with equation \eqref{3.eq} possesses an exponential attractor $\Cal M$ in the sense of Definition \ref{Def2.ea}.
Moreover, the fractal dimension of the global attractor is finite.
\end{teorema}
Indeed, since $\mathscr{K}$ is bounded, the solution $u$ is also automatically bounded in $L^\infty$ (and the same will be true
for the solutions $u_\eb$ of the approximate problems \eqref{1.eqreg} if $\eb>0$ is small enough
 no matter how the {\it regular} interaction function $g$  looks like). So, the term $g(u,x)$ can be treated as a
 perturbation and the proof of Theorem \ref{Th3.intr} repeats word by word the given proof for the particular case
$g(x,u):=-\lambda u$.

\begin{remark}\label{Rem2.grad} The function $g(x,u)$ may even depend explicitly on the gradient $\Nx u$, namely $g=g(x,u,\Nx u)$.
However, in that case we already need to impose some restrictions on the growth of $g$ with respect to $\Nx u$ (since
the obstacle potential controls only the $L^\infty$-norm of a solution and the control of its $W^{1,\infty}$-norm should
be then additionally   obtained). In particular, if $g$ does not grow with respect to $\Nx u$, i.e.,
$$
|g(x,u,\Nx u)|\le Q(|u|)
$$
for some monotone function $Q$ independent of $x$ and $\Nx u$, the proof of Theorem \ref{Th3.intr} still repeats
word by word the case of $g(u)=-\lambda u$. However, our {\bf conjecture} here is that the result remains
true under the standard sub-quadratic growth restriction
$$
|g(x,u,\Nx u)|\le Q(|u|)(1+|\Nx u|^q)
$$
with $q<2$.
\end{remark}

\begin{remark} \label{Rem3.scalar} As we have already pointed out, our method of estimating the Lagrange multipliers is
strongly based on the maximum principle for the leading linear part of equation \eqref{eq}. For this reason, we are
unable in general to extend it to the case of non-scalar diffusion matrix. However, we point out that for some particular convexes
Assumption $\mathscr{L}$ can be verified also for diffusion matrices (say diagonal). This is the case of the convex
$$
\mathscr{K}:=[0,L]^n, \;\;\;L>0.
$$
In this case,  Assumption $\mathscr{L}$ can be easily verified (just multiplying the $k$th equation of \eqref{1.eqreg}
by $\sgn(u^1_k-u_k^2)$ even in the case of diagonal diffusion matrix
\begin{equation}\label{3.eq11}
\partial_t u -a\Delta_{x} u + \partial I_{K}(u) -\lambda u\ni 0,
\end{equation}
with $a=\operatorname{diag}(a_1,\cdots,a_k)$, $a_i>0$ for $i=1,\cdots,n$.
\par
One more generalization can be obtained  replacing the Laplacian $\Dx u$ in equations \eqref{eq}
by the quasi-linear second order differential operator
$$
A(u):=\operatorname{div} (a(|u|^2)\Nx u),\ \ u=(u_1,\cdots,u_n),\ \ |u|^2:=u_1^2+\cdots+u_n^2
$$
with some natural assumptions on the scalar diffusion function $a$. Then, it is not difficult to verify
that the proofs of  Assumption $\mathscr{L}$ given above remain true and the associated solution semigroup possesses
an exponential attractor under the assumptions of Theorem \ref{Th3.intr}.
\end{remark}

\section{Approximations of exponential attractors: the case of a simplex}\label{s4}
The aim of this section is to show that the exponential attractor $\Cal M$ of the singular
problem \eqref{eq} can be approximated by the sequence of exponential attractors $\Cal M_\eb$
of the regular equations \eqref{1.eqreg}. For simplicity, we consider below only the case
in which $\K$ is an $n$-dimensional simplex (the case of an arbitrary bounded convex set with {\it smooth}
boundary can be considered analogously).
To be more precise, the main result of this section is the following Theorem.

\begin{teorema}\label{Th4.main} Let the assumptions of Theorem \ref{Th3.linf} hold. Then, for any
sufficiently small $\eb\ge0$, there exists an exponential attractor $\Cal M_\eb$ for the
 semigroup $S_\eb(t):\Phi_\eb\to\Phi_\eb$
associated with the approximation problem \eqref{1.eqreg} (and the case $\eb=0$ corresponds to the semigroup
$S(t)=S_0(t):\Phi\to\Phi$ associated with the limit singular problem \eqref{eq}). Moreover, the following conditions
hold:
\begin{itemize}
\item{The exponential attractors $\Cal M_\eb$ are uniformly bounded in $C^{2-\nu}$ for any $\nu>0$.}
\item{ $S_\eb(t)\Cal M_\eb\subset \Cal M_\eb$ and the  rate of convergence is uniform with respect to $\eb\to0$:
\begin{equation}\label{4.rate}
\dist_{L^\infty}(S_\eb(t)\Phi_\eb,\Cal M_\eb)\le Ce^{-\alpha t},\ \ \eb\ge0,
\end{equation}
where the positive constants $C$ and $\alpha$ are independent of $\eb$.}
\item{ The fractal dimension  also remains bounded as $\eb\to0$:
\begin{equation}\label{4.dim}
\dim_f(\Cal M_\eb,L^\infty)\le C<\infty,
\end{equation}
where $C$ is independent of $\eb$.}
\item{ The family $\Cal M_\eb$ is H\"older continuous at $\eb=0$:
\begin{equation}\label{4.hol}
\dist^{sym}_{L^\infty}(\Cal M_\eb,\Cal M_0)\le C\eb^\kappa,
\end{equation}
where the positive constants $C$ and $\kappa$ are independent of $\eb$ and $\dist^{sym}$ stands for the symmetric
Hausdorff distance between sets.}
\end{itemize}
\end{teorema}
\begin{proof}
The proof of this Theorem is based on the following abstract result taken
from \cite{FGMZ} (see also \cite{EMZ1}).
\begin{prop}\label{Prop4.abs} Let $\mathscr{H}$ and $\mathscr{H}_1$ be two Banach spaces such that
$\mathscr{H}_1$ is compactly embedded in $\mathscr{H}$ and let, for any $\eb>0$ there exists a bounded closed
set $\Bbb B_\eb\subset \mathscr{H}$, a map $\Bbb S_\eb:\Bbb B_\eb\to\Bbb B_\eb$ and two (nonlinear) projectors
$\Pi_\eb:\Bbb B_\eb\to\Bbb B_0$ and $\Bbb Q_\eb:\Bbb B_0\to\Bbb B_\eb$ such that
the following properties hold:
\par
1) Uniform smoothing property: for every $h_1,h_2\in\Bbb B_\eb$,
\begin{equation}\label{4.sqz}
\|\Bbb S_\eb h_1-\Bbb S_\eb h_2\|_{\mathscr{H}_1}\le L\|h_1-h_2\|_{\mathscr{H}},
\end{equation}
where   $L>0$ is independent of $\eb\ge0$ and $h_1,h_2\in \Bbb B_\eb$.
\par
2) The maps $\Bbb S_\eb$ and $\Bbb S_0$ are close in the following sense:
\begin{equation}\label{4.close}
\begin{cases}
\|\Bbb S_\eb\circ\Bbb Q_\eb h-\Bbb S_0 h\|_{\mathscr{H}}\le C\eb,\ \ \forall h\in\Bbb B_0,\\
\|\Bbb S_0\circ\Pi_\eb h-\Bbb S_\eb h\|_{\mathscr{H}}\le C\eb,\ \ \forall h\in\Bbb B_\eb,
\end{cases}
\end{equation}
where the constant $C$ is independent of $\eb$ and $h$.
\par
Then, the discrete semigroups $\Bbb S_\eb(n)$ generated by the maps $\Bbb S_\eb$ possess
a uniform family of exponential attractors $\Cal M_\eb\subset\Bbb B_\eb$ such that
the following properties hold:
\par
1) Uniform rate of attraction:
\begin{equation}\label{4.drate}
\dist_{\mathscr{H}}(\Bbb S_\eb(n)\Bbb B_\eb,\Cal M_\eb)\le Ce^{-\alpha n}
\end{equation}
for some positive $C$ and $\alpha$ independent of $\eb$.
\par
2) Uniform bounds for the fractal dimension:
\begin{equation}\label{4.ddim}
\dim_f(\Cal M_\eb,\mathscr{H})\le C,
\end{equation}
where $C$ is independent of $\eb$.
\par
3) H\"older continuity at $\eb=0$:
\begin{equation}\label{4.dhol}
\dist^{sym}_{\mathscr{H}}(\Cal M_\eb,\Cal M_0)\le C\eb^\kappa
\end{equation}
for some positive $C$ and $\kappa$ independent of $\eb$.
\end{prop}
We can now prove Theorem \ref{Th4.main} by combining (as we did in the proof of Theorem \ref{Th2.exp})
the Proposition above with the $\ell$ trajectory approach. Namely, we define the spaces $\mathscr{H}$
 and $\mathscr{H}_1$ by \eqref{2.spaces}
and introduce the trajectory
phase spaces
\begin{equation}\label{4.tr}
\Bbb B_\eb:=\{u\in \mathscr{H}_1,\ \ u(t)=S_\eb(t)u_0,\ \ t\in[0,1],\ \ u_0\in\Phi_\eb\}
\end{equation}
and, using the lifting solution operator
$$
\Bbb T_\eb:\Phi_\eb\to\Bbb B_\eb,\ \ (\Bbb T_\eb u_0)(t):=S_\eb(t)u_0,
$$
we lift the solution semigroup $S_\eb(t)$ to the trajectory phase space $\Bbb B_\eb$:
$$
\Bbb S_\eb(t):=\Bbb T_\eb\circ S_\eb(t)\circ \Bbb T_\eb^{-1},\ \ \Bbb S_\eb(t):\Bbb B_\eb\to\Bbb B_\eb
$$
and set, finally, $\Bbb S_\eb:=\Bbb S_\eb(1)$.
\par
Let us verify the assumptions of the abstract proposition for the maps thus defined. To this end, we note
that, thanks to Theorem \ref{Th3.simplex}, we have the uniform estimate \eqref{3.Lest} for the approximations
to the Lagrange multipliers $h_{u_\eb^1}-h_{u^2_\eb}$ associated with two trajectories $u^i_\eb:=S_\eb(t)u_0^i$,
$i=1,2$. Thereofere, arguing exactly as in the proof of Lemma \ref{Lem2.sqz}, we conclude that
\begin{equation}\label{4.sqz1}
\|u^1_\eb-u^2_\eb\|_{L^2(1,2;H^1_0)}+\|\Dt u^1_\eb-\Dt u^2_\eb\|_{L^1(1,2;H^{-s})}\le
L\|u^1_\eb-u^2_\eb\|_{L^2(0,1;L^2)}
\end{equation}
for some positive $L$ which is independent of $\eb$ and $u^i_\eb$, $i=1,2$.
Thus, assumption \eqref{4.sqz} is verified.
\par
In order to verify \eqref{4.close}, we first fix an arbitrary interior point $w_0\in\R^n$ of
$K_0:=K$ and introduce the linear contraction map $\tilde E_\eb: K_\eb\to K$ (where $K_\eb$ is defined by \eqref{3.inv})
by the following expression:
\begin{equation}\label{4.contr}
\tilde E_\eb(u):=\frac1{1+r\eb}(u-w_0)+w_0=\frac1{1+r\eb}u+\frac{r\eb}{1+r\eb}w_0,
\end{equation}
where $r$ is a sufficiently large (but independent of $\eb$) positive constant. Indeed, it is not
difficult to show using the explicit expression for the set $K_\eb$ and for the function $F_\eb$, that
$E_\eb(K_\eb)\subset K_0$ if $r>0$ is large enough.
Finally, we define the map $E_\eb:\Phi_\eb\to\Phi$
as follows
$$
E_\eb(u)(x):=\varphi_\eb(x)\tilde E_\eb(u(x)),
$$
where the cut-off function $\varphi_\eb(x)\in C^1(\R^N)$ is such that $0\le \varphi(x)\le1$, $\varphi_\eb(x)=0$
if $x\in\partial\Omega$ and $\varphi_\eb(x)=0$ if $\dist(x,\partial\Omega)>\eb$. Then, obviously,
\begin{equation}\label{4.closeE}
\|(1-E_\eb)v\|_{L^2}\le C\eb\|v\|_{L^2},\ \ \forall v\in\Phi_\eb
\end{equation}
and for some positive $C$ independent of $\eb$ and $v$.\\
\noindent We now need the following crucial lemma.
\begin{lemma}\label{Lem4.main} Let the assumptions of Theorem \ref{Th4.main} and let $E_\eb$ be defiend by
 \eqref{4.contr}. Then, for every $u_\eb^0\in\Phi_\eb$ and every $u_0\in K_0$, the following estimate holds:
\begin{equation}\label{4.close1}
\|S_\eb(t)u^0_\eb-S_0(t)u_0\|_{L^2}\le C(\eb^{1/2}+\|u_0-u^0_\eb\|_{L^2})e^{Kt},
\end{equation}
where the positive constants $C$ and $K$ are independent of $\eb$, $u_\eb^0\in \Phi_\eb$ and $u_0\in \Phi$.
 \end{lemma}
\begin{proof} Let $\bar{u}(t):=S_0(t)u_0$, $u_\eb(t):=S_\eb(t)u_\eb^0$ and $v_\eb(t):=E_\eb(u_\eb(t))$.
Then, since $v_\eb(t,x)\in K$ for all $t\ge0$, it is an admissible test function for the variational inequality
\eqref{1.vareq}. Therefore,
\begin{equation}\label{4.vareq}
(\Dt \bar{u}(t),\bar{u}(t)-v_\eb(t))+(\Nx \bar{u}(t),\Nx \bar{u}(t)-\Nx v_\eb(t))\le\lambda(\bar{u}(t),\bar{u}(t)-v_\eb(t))
\end{equation}
for almost all $t>0$. Introducing the function $w_\eb(t):=u_\eb(t)-v_\eb(t)=(1-E_\eb)u_\eb(t)$ and \eqref{4.closeE},
we have
\begin{multline}\label{4.vareq1}
(\Dt \bar{u}(t),\bar{u}(t)-u_\eb(t))+(\Nx \bar{u}(t),\Nx \bar{u}(t)-\Nx u_\eb(t))-\lambda(\bar{u}(t),\bar{u}(t)-u_\eb(t))=\\=
(\Dt \bar{u}(t),\bar{u}(t)-v_\eb(t))+(\Nx \bar{u}(t),\Nx \bar{u}(t)-\Nx v_\eb(t))-\lambda(\bar{u}(t),\bar{u}(t)-v_\eb(t))+\\+
(\Dt \bar{u}(t),u_\eb(t)-v_\eb(t))-(\Dx \bar{u}(t), u_\eb(t)-v_\eb(t))-\lambda(\bar{u}(t),u_\eb(t)-v_\eb(t))\le\\
\le -(h_{u_0}(t),w_\eb(t))\le C\eb,
\end{multline}
where we have implicitly used that the Lagrange multiplier $h_{\bar{u}}(t)$ associated with the solution $\bar{u}(t)$
of the singular problem \eqref{1.eq1} is uniformly bounded (see Proposition \ref{Prop1.linf}).
Moreover, multiplying equation \eqref{1.eqreg} by $u_\eb(t)-\bar{u}(t)$, integrating over $x\in\Omega$
and using the monotonicity of $F'_\eb$ and the fact that $F_\eb(\bar{u}(t))\equiv0$, we arrive at
\begin{equation}\label{4.vareq2}
(\Dt u_\eb(t),u_\eb(t)-\bar{u}(t))+(\Nx u_\eb(t),\Nx u_\eb(t)-\Nx \bar{u}(t))-\lambda(u_\eb(t),u_\eb(t)-\bar{u}(t))\le0
\end{equation}
which holds for almost all $t\ge0$. Taking a sum of \eqref{4.vareq1} and \eqref{4.vareq2}, we finally have
\begin{equation}\label{4.fin}
\frac12\frac d{dt}\|u_\eb(t)-\bar{u}(t)\|^2_{L^2}+\|\Nx u_\eb(t)-\Nx \bar{u}(t)\|_{L^2}^2-\lambda\|u_\eb(t)-\bar{u}(t)\|^2_{L^2}\le C\eb
\end{equation}
which, together with the Gronwall inequality, gives \eqref{4.close1} and finishes the proof of the Lemma.
\end{proof}
It is now not difficult to finish the proof of the Theorem.
To this end, we set $\Pi_\eb:\Bbb B_\eb\to\Bbb B_0$
and $\Bbb Q_\eb:\Bbb B_0\to\Bbb B_\eb$ by
$$
\Pi_\eb:= \Bbb T_0\circ E_\eb \circ \Bbb T_\eb^{-1},\ \ \Bbb Q_\eb:=\Bbb T_\eb\circ\Bbb T_0^{-1},
$$
where $\Bbb T_0$ is the lifting operator related to $S_0(t)$.
Then, it is not difficult to see using Lemma \ref{Lem4.main} and estimate \eqref{4.closeE} that the projectors
$\Pi_\eb$ and $\Bbb Q_\eb$ satisfy estimates \eqref{4.close}. Thus, all of the assumptions of the abstract Proposition
\ref{Prop4.abs} are verified and, consequently, the discrete semigroups $\Bbb S_\eb(n)$ acting on the trajectory
phase spaces $\Bbb B_\eb$ possess the uniform family of exponential attractors $\Bbb M_\eb$ which satisfies conditions
\eqref{4.drate}, \eqref{4.ddim} and \eqref{4.dhol} of Proposition \ref{Prop4.abs}.
\par
The rest of the proof can be made exactly as in Theorem \ref{Th2.exp}. Indeed, arguing as in Lemma \ref{Lem2.addreg},
we see that
\begin{equation}\label{2.tr-reg1}
\|u^1_\eb(1)-u^2_\eb(1)\|_{L^2}^2\le (2\lambda+1)\int_0^1\|u^1_\eb(t)-u^2_\eb(t)\|_{L^2}^2\,dt
\end{equation}
for any two solutions $u_\eb^1$ and $u_\eb^2$ of the approximation problems \eqref{1.eqreg}. In addition,
multiplying \eqref{4.fin} by $t$ and integrating over $t\in[0,1]$, we see that
\begin{equation}\label{2.tr-reg2}
\|u_\eb(1)-\bar{u}(1)\|_{L^2}^2\le (2\lambda+1)\int_0^1\|u_\eb(t)-\bar{u}(t)\|_{L^2}^2\,dt.
\end{equation}
Thus, projecting the constructed exponential attractors $\Bbb M_\eb$ back to the physical
phase spaces $\Phi_\eb$ by
$$
\Cal M_\eb^d:=\Bbb M_\eb\big|_{t=1},
$$
we obtain the uniform family of exponential attractors for the discrete semigroups $S_\eb(n):\Phi_\eb\to\Phi_\eb$
which satisfy \eqref{4.rate}, \eqref{4.dim} and \eqref{4.hol} in a weaker space $L^2$ instead of $L^\infty$. In addition,
we have also that
$$
\Cal M_\eb^d\subset S_\eb(1)\Phi_\eb
$$
and, therefore, due to estimate \eqref{3.linf} and to Corollary \ref{Cor1.c2} these attractors are uniformly
bounded in $C^{2-\nu}(\Omega)$ for all $\nu>0$ and the maps $(t,u_0)\to S_\eb(t)u_0$ are uniformly Lipschitz continuous
on $[0,1]\times\Cal M_\eb^d$ (with the Lipschitz constant independent of $\eb$ as well). Thus, the standard formula
$$
\Cal M_\eb:=\cup_{t\in[0,1]}S_\eb(t)\Cal M_\eb^d
$$
gives the desired uniform family of exponential attractors for the semigroups $S_\eb(t):\Phi_\eb\to\Phi_\eb$
with continuous time. This family of exponential attractors
satisfies \eqref{4.rate}, \eqref{4.dim} and \eqref{4.hol} in a weaker topology of $L^2(\Omega)$
instead of $L^\infty(\Omega)$ but, using the fact that $\Cal M_\eb$ are uniformly bounded in $C^{2-\nu}(\Omega)$ together with the interpolation inequality
\eqref{2.int}, we see that \eqref{4.rate}, \eqref{4.dim} and \eqref{4.hol} hold in the initial topology of $L^\infty$
as well. Theorem \ref{Th4.main} is thus completely proved.
\end{proof}

\begin{remark}\label{Rem4.phase} It is not difficult to show that, under the assumptions of Theorem \ref{Th4.main}, the solution
 semigroups $S_\eb(t)$ are actually
defined not only on the space $\Phi_\eb$, but in much lager phase spaces, namely, on $L^2$. Moreover,
arguing analogously to the proof of Theorem \ref{Th3.linf}, one can see that the set $\Phi_\eb$ is an absorbing
set for the semigroup $S_\eb(t)$. This means that, for every bounded $B\subset L^2(\Omega)$ there exists $T=T(\|B\|)$
(independent of $\eb$) such that
$$
S_\eb(t)B\subset\Phi_\eb,\ \ t\ge T.
$$
Thus, the constructed exponential attractors $\Cal M_\eb$ attract (uniformly with respect to $\eb$) the bounded sets
of $L^2(\Omega)$ as well.
\par
Finally, one may even extend the solution semigroup $S(t)$ associated with the obstacle equation \eqref{eq} on
the whole space $L^2$ as well. To this end, we just need to pass to the limit $\eb\to0$ in the solutions
of the approximate problems \eqref{1.eqreg}. It is not difficult to show that the associated
solution $u(t):=S(t)u_0$ will belong to $\Phi$ for every $t>0$ and, in the case
 $u_0\notin\Phi$, it will  have a jump at $t=0$ and the $L^2$-limit
$$
\tilde u_0:= \lim_{t\to0+}u(t)\in\Phi
$$
will exist. Thus, we factually have
$$
u(t)=S(t)u_0=S(t) \Pi_\Phi(u_0),
$$
where $\Pi_\Phi:L^2(\Omega)\to\Phi$ is a non-linear "projector" to the set $\Phi$.
\par
However, in contrast to the semigroup $S(t):\Phi\to\Phi$ which is uniquely defined by the singular equation \eqref{eq},
the above non-linear "projector" $\Pi_\Phi$ essentially depends on the particular choice of the approximation
of the obstacle potential and, therefore, is not canonically defined by the problem \eqref{eq} itself.
For this reason, we have
preferred not to use such projectors in our paper and to define the solution $u(t)$ for the initial data $u(0)\in\Phi$
only.
\end{remark}

\end{document}